\documentclass{article}
\usepackage{amsmath,amssymb,amsthm}
\def\Cmath{\mathbb{C}}
\def\Zmath{\mathbb{Z}}

\newtheorem{lemma}{Lemma}
\newtheorem{proposition}{Proposition}
\newtheorem{theorem}{Theorem}
\newtheorem{corollary}{Corollary}
\newtheorem{definition}{Definition}
\newtheorem{example}{Example}
\newtheorem{remark}{Remark}
\title{Some remarks about the second Leibniz cohomology group for Lie algebras
}
\author{L. Magnin \\
Institut Mathématique de Bourgogne
\thanks{UMR CNRS 5584, Universit\'{e} de Bourgogne, BP 47870, 21078 Dijon Cedex, France.}
\\
\texttt{magnin@u-bourgogne.fr}
}
\date{\today}
\begin{document}
\maketitle
\begin{abstract}
We compare by a very elementary approach the second adjoint and trivial Leibniz
cohomology spaces  of a Lie algebra to the usual ones. Examples are given of coupled cocycles.
Some properties are deduced as to Leibniz deformations.
We also consider the class of Lie algebras for which the Koszul 3-form is zero,
and prove that it contains
all quotients of Borel subalgebras,
or of their nilradicals,
of finite dimensional semisimple Lie algebras.
Finally, a list of Kac-Moody types for indecomposable nilpotent Lie algebras of
dimension $\leq 7$ is given.

\end{abstract}
\section{Introduction.}
Leibniz algebras, along with their Leibniz cohomologies,
were introduced in \cite{loday} as a non antisymmetric version of Lie algebras.
Lie algebras are special Leibniz algebras, and
Pirashvili  introduced \cite{pirashvili} a spectral sequence, that, when applied
to Lie algebras,
measures the  difference between the Lie  algebra cohomology
and the Leibniz cohomology.
Now, Lie algebras have deformations as Leibniz algebras
and those are piloted by the adjoint Leibniz 2-cocycles. In the present paper, we focus
on the second Leibniz cohomology groups
$HL^2(\mathfrak{g},\mathfrak{g}),$
$HL^2(\mathfrak{g},\Cmath)$
for adjoint and trivial
representations of a complex Lie algebra
$\mathfrak{g}$.
We  adopt a very elementary approach,
not resorting to  the Pirashvili sequence, to compare
$HL^2(\mathfrak{g},\mathfrak{g})$ and $HL^2(\mathfrak{g},\Cmath)$
to
$H^2(\mathfrak{g},\mathfrak{g})$ and $H^2(\mathfrak{g},\Cmath)$ respectively.
In both cases,
$HL^2$ appears to be the direct sum of 3 spaces:
$H^2 \oplus ZL^2_0 \oplus \mathcal{C}$ where
$H^2$ is the Lie algebra cohomology group,
$ZL^2_0$ is the space of symmetric
Leibniz-2-cocycles
and $\mathcal{C}$ is a space of \textit{coupled}
Leibniz-2-cocycles the nonzero elements of which have the property that
their symmetric and antisymmetric
parts are not Leibniz cocycles (a similar situation occurs
in second cohomology of current algebras \cite{wagemann}).
Many Lie algebras are uncoupling in the sense that
the space of coupled cocycles is zero.
We give examples in dimensions starting from 5 of
Lie algebras that are not uncoupling.
We deduce some criteria for a Lie algebra to be rigid as Leibniz algebra,
and results as to the number of Leibniz deformations in some cases.
In a first part, we consider the subclass of uncoupling Lie algebras
which have the property that
the Koszul 3-form is zero.
We prove in particular that it contains   all quotients of Borel subalgebras,
or of their nilradicals,
of finite dimensional semisimple Lie algebras.
We also give a list of Kac-Moody types for indecomposable nilpotent Lie algebras of
dimension $\leq 7.$
\section{Leibniz cohomology.}
Recall that a (right) Leibniz algebra is an algebra
$\mathfrak{g}$
with a (non necessarily antisymmetric) bracket,
such that the right adjoint operations $[\cdot,Z]$ are required to be derivations
for any $Z \in \mathfrak{g}.$
In the presence of antisymmetry, that is equivalent to the Jacobi identity, hence any Lie
algebra is a Leibniz algebra.

The Leibniz cohomology
$HL^\bullet(\mathfrak{g},\mathfrak{g})$
of a Leibniz algebra is defined from the complex
$CL^\bullet(\mathfrak{g},\mathfrak{g}) =  \text{Hom }
\left(\mathfrak{g}^{\otimes \bullet}, \mathfrak{g} \right) =
\mathfrak{g} \otimes \left(\mathfrak{g}^*\right)^{\otimes \bullet}$
with the Leibniz-coboundary $\delta$ defined for $\psi \in
CL^n(\mathfrak{g},\mathfrak{g})$ by
\begin{multline*}
(\delta \psi)(X_1,X_2, \cdots , X_{n+1}) =
\\
 [X_1,\psi(X_2,\cdots, X_{n+1})]
 +
 \sum_{i=2}^{n+1} \,(-1)^i [\psi(X_1, \cdots, \hat{X_i}, \cdots, X_{n+1}), X_i]
 \\
 +
 \sum_{1\leqslant i < j \leqslant n+1} (-1)^{j+1}  \,
 \psi(X_1, \cdots, X_{i-1},[X_i,X_j],X_{i+1},\cdots, \hat{X_j}, \cdots, X_{n+1}).
\end{multline*}
(If $\mathfrak{g}$ is a Lie algebra, $\delta$ coincides with the usual coboundary $d$
on $C^\bullet(\mathfrak{g},\mathfrak{g}) =
\mathfrak{g} \otimes \bigwedge^{\bullet} \, \mathfrak{g}^*.$
)

For $\psi \in
CL^1(\mathfrak{g},\mathfrak{g}) =
C^1(\mathfrak{g},\mathfrak{g}) =
\mathfrak{g} \otimes {\mathfrak{g}}^*$
$$(\delta \psi)(X,Y) = [X,\psi(Y)] +[\psi(X),Y] -\psi([X,Y]).$$

For $\psi \in
CL^2(\mathfrak{g},\mathfrak{g}) =
\mathfrak{g} \otimes \left(\mathfrak{g}^*\right)^{\otimes 2},$
\begin{multline*}
(\delta \psi)(X,Y,Z) = [X,\psi(Y,Z)] +[\psi(X,Z),Y] - [\psi(X,Y),Z]
\\
-\psi([X,Y],Z)
+\psi(X,[Y,Z])
+\psi([X,Z],Y).
\end{multline*}

In the same way,
the trivial Leibniz cohomology
$HL^\bullet(\mathfrak{g},\Cmath)$
is defined from the complex
$CL^\bullet(\mathfrak{g},\Cmath) =
\left(\mathfrak{g}^*\right)^{\otimes \bullet}$
with the trivial-Leibniz-coboundary $\delta_\Cmath$ defined for $\psi \in
CL^n(\mathfrak{g},\Cmath)$ by
\begin{multline*}
(\delta_\Cmath \psi)(X_1,X_2, \cdots , X_{n+1}) =
\\
 \sum_{1\leqslant i < j \leqslant n+1} (-1)^{j+1}  \,
 \psi(X_1, \cdots, X_{i-1},[X_i,X_j],X_{i+1},\cdots, \hat{X_j}, \cdots, X_{n+1}).
\end{multline*}
If $\mathfrak{g}$ is a Lie algebra, $\delta_\Cmath$ is
the usual coboundary $d_\Cmath$
on $C^\bullet(\mathfrak{g},\Cmath) =
\bigwedge^{\bullet} \, \mathfrak{g}^*.$

For $\psi \in
CL^1(\mathfrak{g},\Cmath) =
\mathfrak{g}^* \otimes \mathfrak{g},$
$$(\delta_{\Cmath} \psi)(X,Y) = -\psi([X,Y]).$$

For $\psi \in
CL^2(\mathfrak{g},\Cmath) =
\left(\mathfrak{g}^*\right)^{\otimes 2},$
\begin{equation*}
(\delta_\Cmath \psi)(X,Y,Z) =
-\psi([X,Y],Z)
+\psi(X,[Y,Z])
+\psi([X,Z],Y).
\end{equation*}

\section{Some properties of the Koszul map.}

Let $\mathfrak{g}$ be any finite dimensional complex Lie algebra.
Recall that a  symmetric bilinear form
$B \in  S^2 \mathfrak{g}^*$
is invariant, i.e.
$B \in \left(S^2 \mathfrak{g}^*\right) ^{\mathfrak{g}}$
if and only if
$B([Z,X],Y) =-B(X,[Z,Y])  \; \forall X,Y,Z \in
\mathfrak{g}.$
The Koszul map \cite{koszul}
$\mathcal{I} \, : \,
\left(S^2 \mathfrak{g}^*\right) ^{\mathfrak{g}}
\rightarrow
\left(\bigwedge^3 \mathfrak{g}^*\right)^{\mathfrak{g}}
\subset Z^3(\mathfrak{g},\Cmath)
$
is defined by $\mathcal{I}(B)=  I_B,$ with
$I_B(X,Y,Z)=  B([X,Y],Z) \; \forall X,Y,Z \in  \mathfrak{g}.$
\begin{lemma}
\label{lemmep}
Denote $\mathcal{C}^2 \mathfrak{g}=[\frak{g},\frak{g}].$
The projection
$\pi \; : \; \mathfrak{g} \rightarrow \mathfrak{g}/\mathcal{C}^2 \mathfrak{g}$
induces an isomorphism
$$\varpi \; : \; \ker{\mathcal{I}} \rightarrow
S^2  \left( \mathfrak{g}/\mathcal{C}^2 \mathfrak{g} \right)^*
.$$
\end{lemma}
\begin{proof}
For  $B \in \ker{\mathcal{I}},$
define
$\varpi (B) \in S^2 \left(\mathfrak{g}/\mathcal{C}^2 \mathfrak{g}\right)^*$
by $$\varpi(B)(\pi(X),\pi(Y))= B(X,Y),\; \forall X,Y \in \frak{g}.$$
$\varpi(B)$ is well-defined since for $X,Y,U,V \in \frak{g}$
\begin{eqnarray*}
B(X+[U,V],Y)&=& B(X,Y)+B([U,V],Y)\\
&=& B(X,Y)+I_B(U,V,Y)\\
&=& B(X,Y) \text{ (as $I_B=0$)}.
\end{eqnarray*}
The map $\varpi$ is injective since $\varpi(B)=0$ implies $B(X,Y)=0 \; \forall X,Y \in \frak{g}.$
To prove that it is onto, let
$\bar{B} \in S^2  \left( \mathfrak{g}/\mathcal{C}^2 \mathfrak{g} \right)^*,$ and
let
$B_\pi \in S^2 \mathfrak{g}^*$
defined by $B_\pi(X,Y)=\bar{B}(\pi(X),\pi(Y)).$
Then
$B_\pi ([X,Y],Z) = \bar{B}(\pi([X,Y]),\pi(Z))=
\bar{B}(0,\pi(Z)) =0
\linebreak[4]
\; \forall X,Y,Z \in
\frak{g}, $
hence
$B_\pi \in \left(S^2 \mathfrak{g}^*\right) ^{\mathfrak{g}}$
and  $B_\pi \in \ker{\mathcal{I}}.$
Now, $\varpi(B_\pi)=\bar{B}.$
\end{proof}

From lemma
\ref{lemmep},
$\dim
\left(S^2 \mathfrak{g}^*\right) ^{\mathfrak{g}}
= \frac{p(p+1)}{2}
+\dim \text{Im\,}{\mathcal{I}},$
where $p = \dim
 H^1(\mathfrak{g},\Cmath).$
For reductive
$\mathfrak{g},$
$\dim
\left(S^2 \mathfrak{g}^*\right) ^{\mathfrak{g}}
=
\dim H^3(\mathfrak{g},\Cmath)$
(\cite{koszul}).
Note also that the restriction of
$\delta_\Cmath$
to $\left(S^2 \mathfrak{g}^*\right) ^{\mathfrak{g}}$
is $-\mathcal{I}.$

\begin{definition}
$\mathfrak{g}$ is said to be
 $\mathcal{I}$-null
 (resp.
 $\mathcal{I}$-exact)
if
$\mathcal{I}=0$
 (resp.
$\text{Im\,}{\mathcal{I}} \subset B^3(\frak{g},\Cmath))$.
\end{definition}

$\mathfrak{g}$ is $\mathcal{I}$-null if and only
$\mathcal{C}^2 \mathfrak{g} \subset \ker{B}$
$\forall B \in   \left(S^2 \mathfrak{g}^*\right)^{\mathfrak{g}}.$
It is standard that for any
$B \in  \left(S^2 \mathfrak{g}^*\right)^{\mathfrak{g}},$ there exists
$B_1 \in  \left(S^2 \mathfrak{g}^*\right)^{\mathfrak{g}}$
such that
$\ker{(B+B_1)} \subset \mathcal{C}^2 \mathfrak{g},$ hence
$\bigcap_{B \in
\left(S^2 \mathfrak{g}^*\right)^{\mathfrak{g}}}\, \ker{B} \subset \mathcal{C}^2 \mathfrak{g},$
hence
$\mathfrak{g}$ is $\mathcal{I}$-null if and only
$\bigcap_{B \in
\left(S^2 \mathfrak{g}^*\right)^{\mathfrak{g}}}\, \ker{B} = \mathcal{C}^2 \mathfrak{g}.$

\begin{lemma}
(i) Any quotient of a  (not necessarily finite dimensional)
$\mathcal{I}$-null Lie algebra is
$\mathcal{I}$-null.
\\ (ii) any finite direct product of
$\mathcal{I}$-null Lie algebras is
$\mathcal{I}$-null.
\end{lemma}
\begin{proof}
(i)
Let $\frak{g}$ be any
$\mathcal{I}$-null Lie algebra,
 $\frak{h}$ an ideal of
$\frak{g},$
$\bar{\frak{g}}  = \frak{g}/\frak{h},$
\linebreak[4]
$\pi \, : \,
\frak{g}  \rightarrow  \bar{\frak{g}}$ the projection,
and $\bar{B}
\in \left( S^2 {\bar{\mathfrak{g}}}^* \right) ^ {\bar{\frak{g}}}.$
Define $B_\pi
\in S^2 {{\mathfrak{g}}}^*$
by $B_\pi (X,Y) = \bar{B}(\pi(X),\pi(Y)) \, , X,Y \in
\frak{g}.$
Then
$B_\pi ([X,Y],Z) = \bar{B}(\pi([X,Y]),\pi(Z))=
\linebreak[4]
\bar{B}([\pi(X),\pi(Y)],\pi(Z))
= \bar{B}(\pi(X),[\pi(Y),\pi(Z)])
= \bar{B}(\pi(X),\pi([Y,Z]))=
\linebreak[4]
B_\pi (X,[Y,Z]) \,
\forall X,Y,Z \in
\frak{g}, $
hence
$B_\pi \in \left(S^2 \mathfrak{g}^*\right) ^{\mathfrak{g}}$
and
$I_{\bar{B}} \circ (\pi \times \pi \times \pi) =  I_{B_\pi} =0$
since
$\frak{g}$
is $\mathcal{I}$-null.
Hence $I_{\bar{B}}=0.$
\\ (ii)
Let
$\frak{g}  = \frak{g}_1 \times \frak{g}_2$
($ \frak{g}_1 , \frak{g}_2$
 $\mathcal{I}$-null) and
$B\in  \left(S^2 \mathfrak{g}^*\right) ^{\mathfrak{g}}.$
As $B(X_1,[Y_2,Z_2])=
B([X_1,Y_2],Z_2)=
B(0,Z_2)=0\,
\forall X_1 \in  {\frak{g}}_1 ,
Y_2,Z_2 \in  {\frak{g}}_2 ,$
$B$ vanishes on ${\frak{g}}_1 \times {\mathcal{C}}^2{\frak{g}}_2 $
and on $ {\mathcal{C}^2}{\frak{g}}_1 \times  {\frak{g}}_2$ as well,
hence $I_B=0.$
\end{proof}
\begin{lemma}
\label{lemmeborel}
Let $\frak{g}$ be a finite dimensional semi-simple Lie algebra, with Cartan subalgebra
$\frak{h},$
simple root system $S,$ positive roots $\Delta_+,$
and root subspaces ${\frak{g}}^\alpha.$
Let $\frak{k}$ be any non empty subspace of  $\frak{h},$
and $\Gamma \subset
\Delta_+$
such that
$\alpha + \beta \in  \Gamma $
for $\alpha, \beta \in \Gamma, \, \alpha + \beta \in \Delta_+.$
Consider $\frak{u} =\frak{k} \oplus \bigoplus_{\alpha \in \Gamma } {\frak{g}}^\alpha.$
\\ (i) Suppose that $\alpha_{ | \frak{k} } \neq 0 \, \forall \alpha \in \Gamma.$
Then
$\frak{u}$ is $\mathcal{I}$-null.
\\ (ii)
Suppose that $\alpha_{ | \frak{k} } = 0 \, \forall \alpha \in \Gamma \cap S,$
and
$\alpha_{ | \frak{k} } \neq 0 \, \forall \alpha \in \Gamma \setminus S.$
Then
$\frak{u}$ is $\mathcal{I}$-null.
\end{lemma}
\begin{proof}
(i)
Let $\frak{u}_+=  \bigoplus_{\alpha \in \Gamma } {\frak{g}}^\alpha,$
and $X_\alpha$ a root vector in
${\frak{g}}^\alpha$ :
${\frak{g}}^\alpha = \Cmath X_\alpha$ $\forall \alpha \in \Gamma.$
Let  $B\in  \left(S^2 \mathfrak{u}^*\right) ^{\mathfrak{u}}.$
First,
$B(H,X)=0\,  \forall H \in \frak{k},  X
\in \frak{u}_+.$
In fact, for any $\alpha \in \Gamma, $  since
there exists $H_\alpha \in \frak{k} $ such that $\alpha(H_ \alpha) \neq 0,$
$B(H, X_\alpha)=
\frac{1}{\alpha(H_\alpha)} \,
B(H, [H_\alpha, X_\alpha])
=
\frac{1}{\alpha(H_\alpha)} \,
B([H, H_\alpha], X_\alpha])
=
\frac{1}{\alpha(H_\alpha)} \,
B(0, X_\alpha)=0.$
Second, that entails that the restriction of $B$ to
$\frak{u}_+ \times  \frak{u}_+$
is zero, since for any
$\alpha, \beta  \in \Gamma,$
$$B(X_\alpha, X_\beta)= \frac{1}{\alpha(H_\alpha)} B([H_\alpha,X_\alpha], X_\beta)
= \frac{1}{\alpha(H_\alpha)} B(H_\alpha,[X_\alpha, X_\beta])
= 0$$
as $[X_\alpha, X_\beta] \in \frak{u}_+.$
Then
$\frak{u}$ is
$\mathcal{I}$-null.
\\ (ii)
In that case,
$X_\alpha \not \in {\mathcal{C}}^2 \frak{u} \; \forall \alpha \in \Gamma \cap S,$
and $\dim{\left(\frak{u}\left/{\mathcal{C}}^2 \frak{u}\right.\right) }= \dim{\frak{k}}
+\#(\Gamma \cap S).$
For $\frak{u}$ to be
$\mathcal{I}$-null, one has to prove that,
for any  $B\in  \left(S^2 \mathfrak{u}^*\right) ^{\mathfrak{u}}:$
\begin{eqnarray}
\label{1}
B(H,X_\beta)=0 \; \forall H \in \frak{k}, \, \beta \in \Gamma \setminus S ;\\
\label{2}
B(X_\alpha,X_\beta)=0 \; \forall
\alpha \in \Gamma \cap S \, , \,  \beta \in \Gamma \setminus S ;\\
\label{3}
B(X_\beta,X_\gamma)=0 \; \forall
\alpha, \beta \in \Gamma \setminus S.
\end{eqnarray}
(\ref{1}) is proved as in case (i).
To prove (\ref{2}), let $H_\beta \in \frak{k}$ such that $\beta(H_\beta) \neq 0.$
Then
\begin{multline*}
B(X_\alpha,X_\beta)=
\frac{1}{\beta(H_\beta)}B(X_\alpha,[H_\beta,X_\beta])=
\frac{1}{\beta(H_\beta)}B([X_\alpha,H_\beta],X_\beta)= \\
-\frac{1}{\beta(H_\beta)}B(\alpha(H_\beta) X_\alpha,X_\beta)=
-\frac{1}{\beta(H_\beta)}B(0,X_\beta)= 0.
\end{multline*}
As to  (\ref{3}),
\begin{multline*}
B(X_\beta,X_\gamma)=
\frac{1}{\beta(H_\beta)}B([H_\beta,X_\beta],X_\gamma)=
\frac{1}{\beta(H_\beta)}B(H_\beta,[X_\beta,X_\gamma])= 0 \;\text{ from (\ref{1}).}
\end{multline*}

\end{proof}

\begin{example}
\rm
Any Borel subalgebra is
$\mathcal{I}$-null.
\end{example}

\begin{proposition}
\label{prop1}
Let $\mathfrak{g}_2$ be a codimension 1 ideal of
the Lie algebra $\mathfrak{g},$
$(x_1, \cdots , x_N)$ a basis of
$\mathfrak{g}$ with $x_1\not \in \mathfrak{g}_2,$
$x_2, \cdots , x_N \in \mathfrak{g}_2,$
$\pi_2 $
 the corresponding projection onto $\mathfrak{g}_2,$
and $(\omega^1, \cdots ,\omega^N)$ denote the dual basis for
$\mathfrak{g}^*.$
Let $B \in   \left(S^2 \mathfrak{g}^*\right)^{\mathfrak{g}},$
and denote
$B_2 \in   \left(S^2 \mathfrak{g}_2^*\right)^{\mathfrak{g}_2}$
the restriction of $B$ to
$\mathfrak{g}_2 \times \mathfrak{g}_2.$
Then:
\\ (i)
\begin{equation}
I_B=
d(\omega^1 \wedge f) + I_{B_2}\circ (\pi_2\times \pi_2 \times \pi_2).
\end{equation}
where
$f = B(\cdot , x_1) \in \mathfrak{g}^*.$
\\ (ii)
Let   $\gamma \in \bigwedge^2 \mathfrak{g}_2^* \subset \bigwedge^2 \mathfrak{g}^*,$
and denote $d_{\mathfrak{g}_2}$ the coboundary operator
of ${\mathfrak{g}_2}.$  Then
\begin{equation}
\label{drho}
d\gamma =
 \omega^1 \wedge \theta_{x_1} (\gamma)
 +
d_{\mathfrak{g}_2} \gamma
\circ (\pi_2 \times \pi_2 \times \pi_2)
\end{equation}
where $\theta_{x_1}$ stands for the coadjoint action of $x_1$ on the cohomology
of $\mathfrak{g}.$
\\ (iii)
Suppose $I_{B_2} \in B^3(\mathfrak{g}_2,\Cmath),$ and
let  $\gamma \in \bigwedge^2 \mathfrak{g}_2^* \subset \bigwedge^2 \mathfrak{g}^*$ such that
$I_{B_2} =
d_{\mathfrak{g}_2} \gamma.$
Then
$I_{B} \in B^3(\mathfrak{g},\Cmath)$
if and only if
$\omega^1 \wedge \theta_{x_1} (\gamma)
\in B^3(\mathfrak{g},\Cmath).$
In particular, the condition
\begin{equation}
\label{condoneta}
\theta_{x_1} (\gamma) =df
\end{equation}
implies $I_B=d\gamma.$
\end{proposition}

\begin{proof}
(i)
For $X,Y,Z \in \mathfrak{g}$ one has
\begin{multline}
\label{recIB}
B([X,Y],Z) = B([\omega^1(X)x_1+\pi_2(X), \omega^1(Y)x_1+\pi_2(Y)],\omega^1(Z)x_1+\pi_2(Z))
\\ = B\left(\omega^1(X)[x_1,\pi_2(Y)]-\omega^1(Y)
[x_1,\pi_2(X)] +[\pi_2(X),\pi_2(Y)], \right.\\
\left.
\omega^1(Z)x_1+\pi_2(Z)\right)
\\ =
\omega^1(X)\omega^1(Z) B([x_1,\pi_2(Y)],x_1)
-\omega^1(Y)\omega^1(Z) B([x_1,\pi_2(X)],x_1)     \\
+\beta (X,Y,Z) +B\left([\pi_2(X),\pi_2(Y)],\pi_2(Z)\right)
\\
=\beta (X,Y,Z) +B\left([\pi_2(X),\pi_2(Y)],\pi_2(Z)\right)
\end{multline}
where
\begin{multline*}
\beta (X,Y,Z)=
\omega^1(Z) B([\pi_2(X),\pi_2(Y)],x_1)
+\omega^1(X) B([x_1,\pi_2(Y)],\pi_2(Z)) \\
-\omega^1(Y) B([x_1,\pi_2(X)],\pi_2(Z))
\\ =
\omega^1(Z) B([\pi_2(X),\pi_2(Y)],x_1)
+\omega^1(X) B(x_1,[\pi_2(Y),\pi_2(Z)]) \\
-\omega^1(Y) B(x_1,[\pi_2(X),\pi_2(Z)])
.
\end{multline*}
Now
\begin{eqnarray*}
df(X,Y)&=&-B([X,Y],x_1)\\ &=&
-B([\omega^1(X)x_1+\pi_2(X),\omega^1(Y)x_1+\pi_2(Y)],x_1)
\\
&=&
-B(\omega^1(X)[x_1,\pi_2(Y)])-\omega^1(Y)[x_1,\pi_2(X)]
+[\pi_2(X),\pi_2(Y)],x_1)
\\
&=&
-B([\pi_2(X),\pi_2(Y)],x_1),
\end{eqnarray*}
hence
\begin{eqnarray*}
\beta (X,Y,Z)&=&
-(\omega^1(Z) df(X,Y)
+\omega^1(X) df(Y,Z)
-\omega^1(Y) df(X,Z)) \\&=& -(\omega^1 \wedge df)(X,Y,Z).
\end{eqnarray*}
Since $d\omega^1=0,$
 (\ref{recIB}) then reads
\begin{equation}
I_B=
d(\omega^1 \wedge f) + I_{B_2}\circ (\pi_2\times \pi_2 \times \pi_2).
\end{equation}
\\(ii)
One has for any $X,Y,Z \in  \mathfrak{g}$
\begin{multline*}
d\gamma (X,Y,Z) =
d\gamma (\pi_2(X),\pi_2(Y),\pi_2(Z))
+\omega^1(X)d\gamma (x_1,\pi_2(Y),\pi_2(Z))
\\
+\omega^1(Y)d\gamma (\pi_2(X),x_1,\pi_2(Z))
+\omega^1(Z)d\gamma (\pi_2(X),\pi_2(Y),x_1).
\end{multline*}
Now, since $\gamma$ vanishes if one of its arguments is $x_1,$
\begin{eqnarray*}
d\gamma (x_1,\pi_2(Y),\pi_2(Z))
&=&-\gamma ([x_1,\pi_2(Y)],\pi_2(Z)) +\gamma ([x_1,\pi_2(Z)],\pi_2(Y))\\
d\gamma (\pi_2(X),x_1,\pi_2(Z))
&=&-\gamma ([\pi_2(X),x_1],\pi_2(Z)) -\gamma ([x_1,\pi_2(Z)],\pi_2(X))\\
d\gamma (\pi_2(X),\pi_2(Y),x_1)
&=&\gamma ([\pi_2(X),x_1],\pi_2(Y)) -\gamma ([\pi_2(Y),x_1],\pi_2(X)),
\end{eqnarray*}
hence
\begin{multline*}
d\gamma (X,Y,Z) =
d\gamma (\pi_2(X),\pi_2(Y),\pi_2(Z))
+\omega^1(X) \theta_{x_1}\gamma (\pi_2(Y),\pi_2(Z))\\
-\omega^1(Y) \theta_{x_1}\gamma (\pi_2(X),\pi_2(Z))
+\omega^1(Z) \theta_{x_1}\gamma (\pi_2(X),\pi_2(Y))
\\
=d\gamma (\pi_2(X),\pi_2(Y),\pi_2(Z))
+\left(\omega^1 \wedge  \theta_{x_1}\gamma \right)(X,Y,Z)
\end{multline*}
since
$\theta_{x_1}\gamma (\pi_2(U),\pi_2(V))
=\theta_{x_1}\gamma (U,V)$
for all $U,V \in  \mathfrak{g}.$
\\(iii) Results immediately from (i) and (ii).
\end{proof}
\par

\begin{corollary}
\label{thecorr}
Under the hypotheses of Prop. 1, suppose that
$x_1$ commutes with every $x_i$ ($2 \leqslant i \leqslant N$)
except for $x_{i_1}, \cdots, x_{i_r}$  and  that
 $x_{i_1}, \cdots, x_{i_r}$  commute to one another.
Then, if
 $\frak{g}_2$ is $\mathcal{I}$-null,
 $\frak{g}$ is $\mathcal{I}$-null.
\end{corollary}
\begin{proof}
One has to prove that for any invariant bilinear symmetric form $B$
on $\frak{g},$ and any $2\leqslant i,j \leqslant N$,
$B(x_1,[x_i,x_j])=0.$
For $i \neq i_1, \cdots, i_r,$ and any $j\geqslant 2,$
$B(x_1,[x_i,x_j])=
B([x_1,x_i],x_j)=
B(0,x_j)=0.$
For $i,j \in \{i_1, \cdots, i_r\},$
$B(x_1,[x_i,x_j])=
B(x_1,0)=0.$
\end{proof}

\begin{corollary}
Any nilpotent standard filiform Lie algebra or any Heisenberg Lie algebra is
$\mathcal{I}$-null.
\end{corollary}

\par
\begin{corollary}
Any Lie algebra  containing some
$\mathcal{I}$-null
codimension 1 ideal
is $\mathcal{I}$-exact.
\end{corollary}

\begin{corollary}
\label{corollarynilradical}
Any  nilradical $\frak{g}$
of a Borel subalgebra
of a finite-dimensional  semi-simple
Lie algebra is
$\mathcal{I}$-null.
\end{corollary}
\begin{proof}
It is enough to consider the case of a simple Lie algebra, hence of one of the 4 classical types
plus the 5 exceptional ones.
\par Case $A_n$.
Denote $E_{i,j}, 1\leqslant i,j \leqslant n+1$ the canonical basis of
$\frak{gl}(n+1,\Cmath).$
One may suppose that the Borel subalgebra of $A_n$
is comprised of the upper triangular matrices with zero trace, and the Cartan subalgebra $\frak{h}$
is $\bigoplus_{i=1}^{i=n}\; \Cmath H_i$
with
$H_i=E_{i,i}-E_{i+1,i+1}.$
For $n=1$,
 $\mathfrak{g}$ is $\Cmath,$ hence
$\mathcal{I}$-null.
Suppose the result holds true for the nilradical of
the Borel subalgebra of
$A_{n-1}= \frak{sl}(n).$ One has
 $\mathfrak{g} =  \Cmath E_{1,2} \oplus  \cdots  \oplus \Cmath E_{1,n+1} \oplus  \mathfrak{g}^{\prime}_2$
 with $\mathfrak{g}^{\prime}_2$
having basis
$E_{i,j}, 2\leqslant i<j \leqslant n+1$ and being the nilradical of
the Borel subalgebra of
$A_{n-1},$
hence $\mathcal{I}$-null.
 $E_{1,n+1}$ commutes with
 $\mathfrak{g}^{\prime}_2,$ hence, from corr. \ref{thecorr},
 $\Cmath E_{1,n+1} \oplus
 \mathfrak{g}^{\prime}_2$ is $\mathcal{I}$-null.
Now $E_{1,n}$ commutes
with  all members of the basis of
$\Cmath E_{1,n+1} \oplus  \mathfrak{g}^{\prime}_2,$
except for $E_{n,n+1}.$
Hence, from corr. \ref{thecorr},
 $\Cmath E_{1,n} \oplus  \Cmath E_{1,n+1} \oplus  \mathfrak{g}^{\prime}_2$
is $\mathcal{I}$-null.
Now consider
 $\Cmath E_{1,n-1}\oplus \Cmath E_{1,n}\oplus  \Cmath E_{1,n+1} \oplus  \mathfrak{g}^{\prime}_2.$
 $\Cmath E_{1,n-1}$ commutes
with  all members of the basis of
 $ \Cmath E_{1,n}\oplus  \Cmath E_{1,n+1} \oplus  \mathfrak{g}^{\prime}_2$
except for $E_{n-1,n},E_{n-1,n+1}.$ Since those 2 commute, we get in the same way
 that $\Cmath E_{1,n-1}\oplus \Cmath E_{1,n}\oplus  \Cmath E_{1,n+1} \oplus  \mathfrak{g}^{\prime}_2$
is $\mathcal{I}$-null.
The result then follows by induction.
\par Case $D_n$.
We may take $D_n$ as the Lie algebra of matrices
\begin{equation}
\label{Dn}
\begin{pmatrix}
Z_1&Z_2\\Z_3&-{}^tZ_1
\end{pmatrix}
\end{equation}
with $Z_i  \in \frak{gl}(n,\Cmath), Z_2,Z_3$ skew symmetric (see \cite{helgason}, p. 193).
Denote $\tilde{E}_{i,j}=
\left(\begin{smallmatrix}
E_{i,j}&0\\0&-E_{j,i}
\end{smallmatrix}  \right),$
$\tilde{F}_{i,j}=
\left(\begin{smallmatrix}
0&E_{i,j}-E_{j,i}\\0&0
\end{smallmatrix}  \right)$
( $E_{i,j}, 1\leqslant i,j \leqslant n$ the canonical basis of
$\frak{gl}(n,\Cmath)$).
The Cartan subalgebra $\frak{h}$
is $\bigoplus_{i=1}^{i=n}\; \Cmath H_i$
with
$H_i=\tilde{E}_{i,i}$ and the nilradical
of the Borel subalgebra
is
$D_n^+=
\bigoplus_{1\leqslant i < j \leqslant n}\; \Cmath \tilde{E}_{i,j}
\oplus  \bigoplus_{1\leqslant i < j \leqslant n}\; \Cmath \tilde{F}_{i,j}.$
All $\tilde{F}_{i,j}$'s commute to one another, and one has:
\begin{equation}
\label{commFtilde}
[\tilde{E}_{i,j},\tilde{F}_{k,l}] = \delta_{j,k} \tilde{F}_{i,l} - \delta_{j,l}\tilde{F}_{i,k}
\end{equation}
We identify $D_{n-1} $ to a subalgebra of $D_n$ by simply taking the first row and
first column of each block to be zero in (\ref{Dn}).
For $n=2$,
$D_2^+=\Cmath^2$ is $\mathcal{I}$-null.
Suppose the result holds true for
$D_{n-1}^+.$ One has
\begin{equation}
\label{Dn+}
D_n^+=
\Cmath \tilde{E}_{1,2} \oplus
\Cmath \tilde{E}_{1,3} \oplus  \cdots
\Cmath \tilde{E}_{1,n} \oplus
\Cmath \tilde{F}_{1,n} \oplus  \cdots
\oplus \Cmath \tilde{F}_{1,2} \oplus  D_{n-1}^+.
\end{equation}
Start with $\Cmath \tilde{F}_{1,2} \oplus  D_{n-1}^+.$
From
(\ref{commFtilde}),
$\tilde{F}_{1,2}$ commutes with
all
$\tilde{E}_{i,j}\,( 2 \leqslant i <j \leqslant n$) hence with  $D_{n-1}^+.$
Then
$D_{n-1}^+$  is  a codimension 1 ideal of
$\Cmath \tilde{F}_{1,2} \oplus  D_{n-1}^+$
and
$\Cmath \tilde{F}_{1,2} \oplus  D_{n-1}^+$
is $\mathcal{I}$-null from  corr. \ref{thecorr}.
Consider now  $\Cmath \tilde{F}_{1,3} \oplus (\Cmath \tilde{F}_{1,2} \oplus  D_{n-1}^+).$
Again from
(\ref{commFtilde}),
$\tilde{F}_{1,3}$ commutes with
all elements of the basis of
$D_{n-1}^+$ except   $\tilde{E}_{2,3}$ and
$[\tilde{E}_{2,3}, \tilde{F}_{1,3}] = \tilde{F}_{1,2}.$
Then
$\Cmath \tilde{F}_{1,2} \oplus  D_{n-1}^+$
is  a codimension 1 ideal of
$\Cmath \tilde{F}_{1,3} \oplus (\Cmath \tilde{F}_{1,2} \oplus  D_{n-1}^+),$
and  the latter is $\mathcal{I}$-null.
Suppose that
$\Cmath \tilde{F}_{1,s-1} \oplus \cdots  \oplus \Cmath \tilde{F}_{1,2} \oplus  D_{n-1}^+$
is a
codimension 1 ideal of
$\Cmath \tilde{F}_{1,s} \oplus (\Cmath \tilde{F}_{1,s-1} \oplus \cdots  \oplus \Cmath \tilde{F}_{1,2} \oplus  D_{n-1}^+),$
 and that the latter is $\mathcal{I}$-null.
Consider
$\Cmath \tilde{F}_{1,s+1} \oplus (\Cmath \tilde{F}_{1,s} \oplus \cdots  \oplus \Cmath \tilde{F}_{1,2} \oplus  D_{n-1}^+).$
From
(\ref{commFtilde}),
for $2 \leqslant i < j \leqslant n,$
$[\tilde{E}_{i,j},\tilde{F}_{1,s+1}]= \delta_{j,s+1} \tilde{F}_{1,i}$
is nonzero only for $i=2, \cdots, s,$ and $j=s+1,$ and it is then equal to
$\tilde{F}_{1,i}.$
Then first
$\Cmath \tilde{F}_{1,s} \oplus (\Cmath \tilde{F}_{1,s-1} \oplus \cdots  \oplus \Cmath \tilde{F}_{1,2} \oplus  D_{n-1}^+)$
is a codimension 1 ideal of
$\Cmath \tilde{F}_{1,s+1} \oplus (\Cmath \tilde{F}_{1,s} \oplus \cdots  \oplus \Cmath \tilde{F}_{1,2} \oplus  D_{n-1}^+).$
Second,
 the latter
is $\mathcal{I}$-null from  corr. \ref{thecorr}.
By induction the above property holds for $s=n.$
Consider now
$\Cmath \tilde{E}_{1,n} \oplus (
\Cmath \tilde{F}_{1,n} \oplus  \cdots
\oplus \Cmath \tilde{F}_{1,2} \oplus  D_{n-1}^+).$
One has for
 $2 \leqslant i <j \leqslant n,$
$[\tilde{E}_{1,n},\tilde{E}_{i,j}]=0,$
$[\tilde{E}_{1,n},\tilde{F}_{i,j}]=   -\delta_{n,j} \tilde{F}_{1,i},$
$[\tilde{E}_{1,n},\tilde{F}_{1,j}]= 0.$
Hence
$\Cmath \tilde{F}_{1,n} \oplus  \cdots
\oplus \Cmath \tilde{F}_{1,2} \oplus  D_{n-1}^+$ is an ideal of
$\Cmath \tilde{E}_{1,n} \oplus (
\Cmath \tilde{F}_{1,n} \oplus  \cdots
\oplus \Cmath \tilde{F}_{1,2} \oplus  D_{n-1}^+)$  and the latter is
$\mathcal{I}$-null.
For $2 \leqslant i <j \leqslant n,$
$1 \leqslant k \leqslant n-2,$
\begin{eqnarray*}
[\tilde{E}_{1,n-k},\tilde{E}_{i,j}]&=& \delta_{n-k,i} \tilde{E}_{1,j}, \cr
[\tilde{E}_{1,n-k},\tilde{F}_{i,j}]&=& \delta_{n-k,i} \tilde{F}_{1,j} - \delta_{n-k,j} \tilde{F}_{1,i},\cr
 [  {\tilde{E}}_{1,n-k} , \tilde{E}_{1,n} ]& =&   \delta_{n-k,1} \tilde{E}_{1,n},\cr
[\tilde{E}_{1,n-k},\tilde{F}_{1,j}]&=& \delta_{n-k,1} \tilde{F}_{1,j}.
\end{eqnarray*}
$[\tilde{E}_{1,n-1},\tilde{E}_{i,j}]$
is nonzero only for $i=n-1,j=n$ and yields then
$\tilde{E}_{1,n};$
\linebreak[4]
$[\tilde{E}_{1,n-1},\tilde{F}_{i,j}]$
is nonzero only for ($i=n-1,j=n$)
or for ($i < j=n-1$)
and yields respectively
$ \tilde{F}_{1,n},$ or
$- \tilde{F}_{1,i}.$
 $[  {\tilde{E}}_{1,n-1} , \tilde{E}_{1,n} ]$ and
$[\tilde{E}_{1,n-1},\tilde{F}_{1,j}]$ are zero for $n\geqslant 3.$
Hence, first
$\Cmath \tilde{E}_{1,n} \oplus
\Cmath \tilde{F}_{1,n} \oplus  \cdots
\oplus \Cmath \tilde{F}_{1,2} \oplus  D_{n-1}^+,$
is a codimension 1 ideal of
$\Cmath \tilde{E}_{1,n-1} \oplus (
\Cmath \tilde{E}_{1,n} \oplus
\Cmath \tilde{F}_{1,n} \oplus  \cdots
\oplus \Cmath \tilde{F}_{1,2} \oplus  D_{n-1}^+),$
and second ,
 the latter is $\mathcal{I}$-null,
since
$ \tilde{E}_{n-1,n},$
commutes with
$ \tilde{F}_{i,n-1},$
$ \tilde{F}_{n-1,n}.$
Suppose that
$\Cmath \tilde{E}_{1,n-k+1} \oplus \cdots
\oplus \Cmath \tilde{E}_{1,n}
\oplus \Cmath \tilde{F}_{1,n} \oplus  \cdots
\oplus \Cmath \tilde{F}_{1,2} \oplus  D_{n-1}^+$
is a
codimension 1 ideal of
$\Cmath \tilde{E}_{1,n-k} \oplus (
\Cmath \tilde{E}_{1,n-k+1} \oplus \cdots
\oplus \Cmath \tilde{E}_{1,n}
\oplus \Cmath \tilde{F}_{1,n} \oplus  \cdots
\oplus \Cmath \tilde{F}_{1,2} \oplus  D_{n-1}^+)$
and that the latter is
 $\mathcal{I}$-null.
Consider
$\Cmath \tilde{E}_{1,n-k-1} \oplus (
\Cmath \tilde{E}_{1,n-k} \oplus \cdots
\oplus \Cmath \tilde{E}_{1,n}
\oplus \Cmath \tilde{F}_{1,n} \oplus  \cdots
\oplus \Cmath \tilde{F}_{1,2} \oplus  D_{n-1}^+).$
$[\tilde{E}_{1,n-k-1},\tilde{E}_{i,j}]$
is nonzero only for $i=n-k-1$ and yields then
$\tilde{E}_{1,j};$
$[\tilde{E}_{1,n-k-1},\tilde{F}_{i,j}]= \delta_{n-k-1,i} \tilde{F}_{1,j} - \delta_{n-k-1,j} \tilde{F}_{1,i}$
is nonzero only for $i=n-k-1$ or $j=n-k-1$ and yields resp.
$\tilde{F}_{1,j}$ or
$-\tilde{F}_{1,i}.$
Hence
$\Cmath \tilde{E}_{1,n-k} \oplus \cdots
\oplus \Cmath \tilde{E}_{1,n}
\oplus \Cmath \tilde{F}_{1,n} \oplus  \cdots
\oplus \Cmath \tilde{F}_{1,2} \oplus  D_{n-1}^+$
is an ideal of
$\Cmath \tilde{E}_{1,n-k-1} \oplus (
\Cmath \tilde{E}_{1,n-k} \oplus \cdots
\oplus \Cmath \tilde{E}_{1,n}
\oplus \Cmath \tilde{F}_{1,n} \oplus  \cdots
\oplus \Cmath \tilde{F}_{1,2} \oplus  D_{n-1}^+).$
The latter
 is $\mathcal{I}$-null since
$\tilde{E}_{n-k-1,j}$ commutes with both
$\tilde{F}_{n-k-1,j^{\prime}},$
$\tilde{F}_{i,n-k-1}.$
The result follows by induction.
\par Case $B_n$.
We may take $B_n$ ($n\geqslant 2$) as the Lie algebra of matrices
\begin{equation}
\left(
\begin{array}{c|c|c}
0& u&v\\
\hline
-{}^tv&Z_1&Z_2\\
\hline
-{}^tu&Z_3&-{}^tZ_1
\end{array}
\right)
\end{equation}
with $u,v$ complex $(1\times n)$-matrices,
$Z_i  \in \frak{gl}(n,\Cmath), Z_2,Z_3$ skew symmetric, i.e.
\begin{equation}
\left(
\begin{array}{c|ccc}
0& u&&v\\
\hline
-{}^tv&& &\\
&&\boxed{A} &\\
-{}^tu&&
\end{array}
\right)
\end{equation}
with $A \in D_n.$
The Cartan subalgebra is the same as for $D_n.$
Then $B_n^+$  consists of the matrices
\begin{equation}
\left(
\begin{array}{c|ccc}
0& 0&&v\\
\hline
-{}^tv&& &\\
&&\boxed{A} &\\
0&&
\end{array}
\right)
\end{equation}
with $v$ complex $(1\times n)$-matrix and $A \in D_n^+.$
We identify $A \in D_n^+$ to the matrix
\begin{equation*}
\left(
\begin{array}{c|ccc}
0& 0&&0\\
\hline
0&& &\\
&&\boxed{A} &\\
0&&
\end{array}
\right)       \in B_n^+.
\end{equation*}
For $1 \leqslant q \leqslant n, $ let
${v_q}$ the $1\times n$-matrix $E_{1,q}=(0, \cdots,1,\cdots,0)$ ($1$ in $q^{\text{th}}$ position), and
\begin{equation*}
\tilde{v}_q
=
\left(
\begin{array}{c|ccc}
0& 0&&v_q\\
\hline
-{}^tv_q&& &\\
&&\boxed{0} &\\
0&&
\end{array}
\right)
\end{equation*}
Hence $B_n^+ = \left( \bigoplus_{q=1}^{n}
\, \Cmath \tilde{v}_q \right) \oplus D_n^+.$
One has for
$1 \leqslant q \leqslant n, \; 1 \leqslant i<j \leqslant n $
\begin{eqnarray*}
[\tilde{v}_q ,\tilde{E}_{i,j}]&=&   -\delta_{q,j} \tilde{v}_i
\cr
[\tilde{v}_q ,\tilde{F}_{i,j}]&=&   0
\end{eqnarray*}
and for $1 \leqslant s<q \leqslant n$
\begin{equation}
[\tilde{v}_q ,\tilde{v}_{s}]=   \tilde{F}_{s,q} .
\end{equation}
Consider
$\Cmath \tilde{v}_1 \oplus D_n^+.$
As
$\tilde{v}_1$ commutes with
$\tilde{E}_{i,j}, \tilde{F}_{i,j},$
$D_n^+$ is an ideal of
$\Cmath \tilde{v}_1 \oplus D_n^+$
and the latter is  $\mathcal{I}$-null.
 Suppose
$\Cmath \tilde{v}_{s-1} \oplus \cdots \oplus \Cmath \tilde{v}_1 \oplus D_n^+$
is an ideal of $\Cmath \tilde{v}_{s} \oplus \left( \Cmath \tilde{v}_{s-1} \oplus \cdots \oplus \Cmath \tilde{v}_1 \oplus D_n^+\right)$
and the latter is  $\mathcal{I}$-null.
Consider
$\Cmath \tilde{v}_{s+1} \oplus \left( \Cmath \tilde{v}_{s} \oplus  \Cmath \tilde{v}_{s-1} \oplus \cdots \oplus \Cmath \tilde{v}_1 \oplus D_n^+\right).$
$[\tilde{v}_{s+1} ,\tilde{E}_{i,j}]=   -\delta_{s+1,j} \tilde{v}_i$ hence
$\tilde{v}_{s+1}$  commutes to all $\tilde{E}_{i,j}$s except for
$\tilde{E}_{i,s+1}$ ($i< s+1$) and then yields
$-\tilde{v}_i.$
For $t <s+1,$
$[\tilde{v}_{s+1} ,\tilde{v}_{t}]=   \tilde{F}_{t,s+1} .$
Hence $\Cmath \tilde{v}_{s} \oplus  \Cmath \tilde{v}_{s-1} \oplus \cdots \oplus \Cmath \tilde{v}_1 \oplus D_n^+$
is an ideal of
$\Cmath \tilde{v}_{s+1} \oplus \left( \Cmath \tilde{v}_{s} \oplus  \Cmath \tilde{v}_{s-1} \oplus \cdots \oplus \Cmath \tilde{v}_1 \oplus D_n^+\right).$
The latter is  $\mathcal{I}$-null
since the
$\tilde{E}_{i,s+1}$ ($i< s+1$) commute to one another.
 Hence by induction the property holds for $s=n$ and $B_n^+$
is  $\mathcal{I}$-null.
\par Case $C_n$.
This case is pretty similar to the case $D_n.$
We may take $C_n$ as the Lie algebra of matrices
\begin{equation}
\label{Cn}
\begin{pmatrix}
Z_1&Z_2\\Z_3&-{}^tZ_1
\end{pmatrix}
\end{equation}
with $Z_i  \in \frak{gl}(n,\Cmath), Z_2,Z_3$ symmetric.
 $\tilde{E}_{i,j}$ and the Cartan subalgebra are identical to those of $D_n.$
We denote for $1 \leqslant i,j \leqslant n:$
$\hat{F}_{i,j}=
\left(\begin{smallmatrix}
0&E_{i,j}+E_{j,i}\\0&0
\end{smallmatrix}  \right).$ Then
\begin{equation}
\label{Cn+}
C_n^+=
\bigoplus_{1\leqslant i < j \leqslant n}\; \Cmath \tilde{E}_{i,j}
\oplus   \bigoplus_{1\leqslant k \leqslant l \leqslant n}\; \Cmath \hat{F}_{k,l}.
\end{equation}
All $\hat{F}_{k,l}$'s commute to one another, and one has:
\begin{equation}
\label{hatF}
[\tilde{E}_{i,j},\hat{F}_{k,l}] = \delta_{j,k} \hat{F}_{i,l} + \delta_{j,l}\hat{F}_{i,k}
\end{equation}
The case is step by step analog to the case of $D_n$ with
(\ref{hatF}) instead of (\ref{commFtilde})
and (\ref{Cn+}) instead of (\ref{Dn+}).
\par Case $G_2$.
The commutation relations for $G_2$  appear in \cite{fulton-harris}, p. 346.
$G_2^+$ is 6-dimensional with comm. rel.
$[x_1,x_2]=x_3; \,
[x_1,x_3]=x_4;   \,
[x_1,x_4]=-3x_5;   \,
[x_2,x_5]=-x_6;   \,
[x_3,x_4]=-3x_6.$
$G_2^+$ has the same adjoint cohomology $(1,4,7,8,7,5,2)$ as,
and is  isomorphic to,
$\mathfrak{g}_{6,18},$ which is
 $\mathcal{I}$-null.
\par Case $F_4$.
$F_4^+$ has  24 positive roots,
and root vectors $x_i$ ($1 \leqslant i \leqslant 24$).
From the root pattern, one gets with some calculations the
commutation relations of $F_4^+:$
$ [x_1,x_2]=x_{5};$
$ [x_1,x_{13}]=x_{14};$
$ [x_1,x_{15}]= - x_{6};$
$ [x_1,x_{16}]= - x_{7};$
$ [x_1,x_{17}]= - x_{23};$
$ [x_1,x_{18}]=x_{19};$
$ [x_1,x_{24}]=x_{22};$
$ [x_2,x_3]=x_{15};$
$ [x_2,x_7]=x_{8};$
$ [x_2,x_{12}]=x_{13};$
$ [x_2,x_{19}]=x_{20};$
$ [x_2,x_{21}]=x_{24};$
$ [x_2,x_{23}]=x_{9};$
$ [x_3,x_4]=x_{21};$
$ [x_3,x_5]=x_{6};$
$ [x_3,x_6]=x_{7};$
$ [x_3,x_9]=x_{10};$
$ [x_3,x_{11}]=x_{12};$
$ [x_3,x_{15}]=x_{16};$
$ [x_3,x_{20}]= - 2 x_{11};$
$ [x_3,x_{22}]=
\frac{1}{2}
x_{23};$
$ [x_3,x_{24}]= -
\frac{1}{2}
x_{17};$
$ [x_4,x_6]=x_{22};$
$ [x_4,x_7]=x_{23};$
$ [x_4,x_8]=x_{9};$
$ [x_4,x_9]= - x_{20};$
$ [x_4,x_{10}]=x_{11};$
$ [x_4,x_{15}]= - x_{24};$
$ [x_4,x_{16}]=x_{17};$
$ [x_4,x_{17}]=x_{18};$
$ [x_4,x_{23}]= - x_{19};$
$ [x_5,x_{12}]=x_{14};$
$ [x_5,x_{16}]=x_{8};$
$ [x_5,x_{17}]=x_{9};$
$ [x_5,x_{18}]= - x_{20};$
$ [x_5,x_{21}]=x_{22};$
$ [x_6,x_{11}]= - x_{14};$
$ [x_6,x_{15}]= - x_{8};$
$ [x_6,x_{17}]=x_{10};$
$ [x_6,x_{18}]=2 x_{11};$
$ [x_6,x_{21}]=
\frac{1}{2}
x_{23};$
$ [x_6,x_{24}]=
\frac{1}{2}
x_{9};$
$ [x_7,x_{18}]=2 x_{12};$
$ [x_7,x_{20}]= - 2 x_{14};$
$ [x_7,x_{24}]=x_{10};$
$ [x_8,x_{18}]=2 x_{13};$
$ [x_8,x_{19}]=2 x_{14};$
$ [x_8,x_{21}]= - x_{10};$
$ [x_9,x_{17}]= - 2 x_{13};$
$ [x_9,x_{21}]= - x_{11};$
$ [x_9,x_{23}]=2 x_{14};$
$ [x_{10},x_{21}]= - x_{12};$
$ [x_{10},x_{22}]= - x_{14};$
$ [x_{10},x_{24}]= - x_{13};$
$ [x_{11},x_{15}]= - x_{13};$
$ [x_{15},x_{19}]=2 x_{11};$
$ [x_{15},x_{21}]=
\frac{1}{2}
x_{17};$
$ [x_{15},x_{22}]=
\frac{1}{2}
x_{9};$
$ [x_{15},x_{23}]= - x_{10};$
$ [x_{16},x_{19}]=2 x_{12};$
$ [x_{16},x_{20}]=2 x_{13};$
$ [x_{16},x_{22}]=x_{10};$
$ [x_{17},x_{22}]=x_{11};$
$ [x_{17},x_{23}]=2 x_{12};$
$ [x_{21},x_{22}]=
\frac{1}{2}
x_{19};$
$ [x_{21},x_{24}]=
\frac{1}{2}
x_{18};$
$ [x_{22},x_{24}]=
 -
\frac{1}{2}
 x_{20};$
$ [x_{23},x_{24}]=x_{11}.$
Then one verifies
(see \cite{companionarchive})
that $F_4^+$
 is $\mathcal{I}$-null.
\par Case $E_6$.
In the case of $E_6^+$ the set $\Delta_+$ of positive roots
(associated to
the set $S$ of simple roots) has cardinality
36  (\cite{fulton-harris}, p. 333):
\begin{multline*}
\Delta_+ =
\{\varepsilon_i+\varepsilon_j; 1\leqslant i <j \leqslant5\} \cup
\{\varepsilon_i-\varepsilon_j; 1\leqslant j <i \leqslant5\}
\\
\cup
\{\frac{1}{2}(
\pm\varepsilon_1\pm\varepsilon_2
\pm\varepsilon_3\pm\varepsilon_4
\pm\varepsilon_5+\sqrt{3}\, \varepsilon_6);
\text{\# minus signs even}
\}.
\end{multline*}
Instead of computing commutation relations, we will use the following
property of $\Delta_+.$
\\  ($\mathcal{P}$):
for $\alpha, \beta, \gamma \in \Delta_+,$
if
$\alpha+\beta \in \Delta_+$
and $\alpha+\gamma \in \Delta_+$
then
$\beta+\gamma \not\in \Delta_+.$
\\
Now, introduce some Chevalley basis (\cite{helgason}, p. 19 ex. 7)
of
$E_6^+:$
$(X_\alpha)_{\alpha \in \Delta^+}.$
One has
$$[X_\alpha, X_\beta] =N_{\alpha,\beta} X_{\alpha+\beta}  \; \forall \alpha, \beta \in
\Delta_+$$
$$N_{\alpha,\beta}=0 \text{ if } \alpha+\beta \not \in \Delta_+,
N_{\alpha,\beta}\in \Zmath \setminus \{0\}  \text{ if } \alpha+\beta  \in \Delta_+.$$
Define inductively a sequence
$\frak{g}_1 \subset
\frak{g}_2 \subset \cdots
\subset
\frak{g}_{36} = E_6^+$
of
 $\mathcal{I}$-null subalgebras, each of which a codimension 1 ideal of the following,
 as follows.
 Start with
$\frak{g}_1 = \Cmath X_{\delta_1},$
$\delta_1 \in \Delta_+$ of maximum height.
Suppose
$\frak{g}_i$ defined. Then take
$\frak{g}_{i+1} = \Cmath X_{\delta_{i+1}} \oplus \frak{g}_i$
with $\delta_{i+1} \in \Delta_+ \setminus \{\delta_1, \cdots ,\delta_i\}$
of maximum height.
Clearly,
$\frak{g}_i$ is a 1-dimensional ideal of $\frak{g}_{i+1}.$
To prove that it is
 $\mathcal{I}$-null
 we only have to check that, for $1 \leqslant s,t \leqslant i,$
 if
 $\delta_{i+1}+\delta_s \in \Delta_+$
 and
 $\delta_{i+1}+\delta_t \in \Delta_+$
 then $\delta_s +\delta_t \not \in \Delta_+.$
 That holds true because of property  ($\mathcal{P}$).
\par Case $E_7$.
In the case of $E_7^+$ the set $\Delta_+$ of positive roots
(associated to
the set $S$ of simple roots) has cardinality
63  (\cite{fulton-harris}, p. 333):
\begin{multline*}
\Delta_+ =
\{\varepsilon_i+\varepsilon_j; 1\leqslant i <j \leqslant6\} \cup
\{\varepsilon_i-\varepsilon_j; 1\leqslant j <i \leqslant6\}
\cup
\{\sqrt{2}\varepsilon_7
\}
\\
\cup
\{\frac{1}{2}(
\pm\varepsilon_1\pm\varepsilon_2
\pm\varepsilon_3\pm\varepsilon_4
\pm\varepsilon_5
\pm\varepsilon_6
+\sqrt{2}\, \varepsilon_7);
\text{\# minus signs odd}
\}.
\end{multline*}
 Property ($\mathcal{P}$) holds true for $E_7^+$
(see \cite{companionarchive}).
 Hence the conclusion follows as in the case of $E_6^+.$
\par Case $E_8$.
In the case of $E_8^+$ the set $\Delta_+$ of positive roots
has cardinality
120 (\cite{fulton-harris}, p. 333):
\begin{multline*}
\Delta_+ =
\{\varepsilon_i+\varepsilon_j; 1\leqslant i <j \leqslant8\} \cup
\{\varepsilon_i-\varepsilon_j; 1\leqslant j <i \leqslant8\}
\\
\cup
\{\frac{1}{2}(
\pm\varepsilon_1\pm\varepsilon_2
\pm\varepsilon_3\pm\varepsilon_4
\pm\varepsilon_5
\pm\varepsilon_6
\pm\varepsilon_7
+ \varepsilon_8);
\text{\# minus signs even}
\}.
\end{multline*}
 Property ($\mathcal{P}$) holds true for $E_8^+$
(see \cite{companionarchive}).
 Hence the conclusion follows as in the case of $E_6^+.$
\end{proof}
\begin{remark}
\rm
 Property ($\mathcal{P}$) holds true for $A_n^+,$ hence we could have used it.
 However, it does not hold true for $F_4^+.$
 One has for example in the above commutation relations of $F_4^+$
 (with root vectors)
 $[x_3,x_4] \neq 0, [x_3,x_9]\neq 0, $ yet $[x_4,x_9] \neq 0.$
\end{remark}

\begin{corollary}
Suppose that  the Lie algebra $\mathfrak{g}$ is such that
$\dim{\text{Im\,}{\mathcal{I}}}=0 \text{ or }1. $
Let $\tau \in \text{Der } \frak{g}$ such that
$\tau x_k \in
\mathcal{C}^2 \mathfrak{g}$ $\forall k\geqslant 2$
where $(x_1, \cdots, x_N)$ is some basis of
$\frak{g}.$ Denote $\tilde{\frak{g}}_{\tau} = \Cmath \, \tau \oplus \frak{g}$
the Lie algebra obtained by adjoining the derivation $\tau$ to $\frak{g},$
and by $\tilde{\mathcal{I}}$ the Koszul map of
$\tilde{\frak{g}}_{\tau}.$
Then  $\dim{\text{Im\,}{\tilde{\mathcal{I}}}}=0 \text{ or }1. $
\end{corollary}
\begin{proof}
Let $B \in \left( S^2 {\tilde{\frak{g}}}_\tau {}^*\right)^{\tilde{\frak{g}}_\tau}$
such that $I_B  \neq  0.$
One has
$I_B= \omega^{\tau}\wedge df_\tau +I_{B_2} \circ (\pi_2 \times \pi_2 \times \pi_2)$
where $(\tau, x_1, \cdots, x_N)$ is the basis of
$\tilde{\frak{g}}_{\tau},$
$(\omega^\tau, \omega^1, \cdots, \omega^N)$  the dual basis, $B_2$ the restriction
of $B$ to $\frak{g}$, $f_\tau = B(\tau, \cdot \, )$ and $\pi_2$ the projection on $\frak{g}.$
One may suppose
$\dim{\text{Im\,}{\mathcal{I}}}=1.$ Let
$C \in \left( S^2 {{\frak{g}}} {}^*\right)^{{\frak{g}}}$
with
$I_C \neq 0.$
There exists $\lambda \in \Cmath$ such that
$I_{B_2} =\lambda I_C.$ Now, for $X,Y \in \frak{g},$
\begin{multline*}
df_\tau(X,Y) = B(\tau,[X,Y]) \\
=
\omega^1(X)B(\tau,[x_1,\pi_3(Y)])
-\omega^1(Y)B(\tau,[x_1,\pi_3(X)])
+B(\tau,[\pi_3(X),\pi_3(Y)])
\\=
-\omega^1(X)B(\tau \pi_3(Y),x_1)
+\omega^1(Y)B(\tau \pi_3(X),x_1)
+B(\tau \pi_3(X),\pi_3(Y))
\end{multline*}
where $\pi_3$ is the projection on $\text{vect}(x_2, \cdots, x_N).$
Now, $B(\tau \pi_3(Y),x_1)=$
\linebreak[4]
$\lambda C(\tau \pi_3(Y),x_1),$
$B(\tau \pi_3(X),x_1)
=\lambda C(\tau \pi_3(X),x_1),$
$B(\tau \pi_3(X),\pi_3(Y))
=$
\linebreak[4]
$\lambda C(\tau \pi_3(X),\pi_3(Y))$
since
$\tau \pi_3(X) ,\tau \pi_3(Y) \in \mathcal{C}^2\frak{g}.$
The result follows.
\end{proof}
\begin{example}
\rm
The nilpotent Lie algebra
$\mathfrak{g}_{7,2.4},$ has commutation relations
$[x_1,x_2]=x_3,$
$[x_1,x_3]=x_4,$
$[x_1,x_4]=x_5,$
$[x_1,x_5]=x_6,$
$[x_2,x_5]=-x_7,$
$[x_3,x_4]=x_7.$
The  elements of $\text{Der } \mathfrak{g}_{7,2.4} (\text{mod } \text{ad} \mathfrak{g}_{7,2.4})$ are
\begin{equation}
\tau=
\begin{pmatrix}
\xi^1_1&0&0&0&0&0&0\\
\xi^2_1&\xi^2_2&0&0&0&0&0\\
0&0&\xi^1_1+\xi^2_2&0&0&0&0\\
0&0&0&2\xi^1_1+\xi^2_2&0&0&0\\
0&\xi^5_2&0&0&3\xi^1_1+\xi^2_2&0&0\\
\xi^6_1&\xi^6_2&\xi^5_2&0&0&4\xi^1_1+\xi^2_2&0\\
\xi^7_1&0&0&0&0&-\xi^1_2&3\xi^1_1+2\xi^2_2
\end{pmatrix}
\end{equation}
$\tau$ is nilpotent if $\xi^1_1=\xi^2_2=0.$
Denote the nilpotent $\tau$ by
$(\xi^2_1;\xi^5_2;\xi^6_1,\xi^6_2;\xi^7_1).$
Now, projectively equivalent derivations $\tau, \tau^{\prime}$ (see \cite{ART})
give isomorphic
$\tilde{\frak{g}}_{\tau},$
$\tilde{\frak{g}}_{\tau^{\prime}}.$
By reduction using projective equivalence, we are reduced to the following cases:
Case 1. $\xi^2_1\neq 0 : (1;\varepsilon;0,\eta;0);$
Case 2. $\xi^2_1= 0 : (0;\varepsilon;0,\eta;\lambda);$
where $\varepsilon, \eta, \lambda= 0, 1.$
In both cases
$\tilde{\frak{g}}_{\tau}$
is $\mathcal{I}$-null,
except when $\tau=0$ in case 2 where
$\tilde{\frak{g}}_{\tau}$
is the direct product
$\Cmath \times \mathfrak{g}_{7,2.4}$ which is quadratic
(i.e. having a nondegenerate invariant bilinear form).  Hence any indecomposable
8-dimensional nilpotent Lie algebra containing a subalgebra isomorphic to
$\mathfrak{g}_{7,2.4}$
is $\mathcal{I}$-null, though
$\mathfrak{g}_{7,2.4}$ is quadratic.
That is in line with the fact  that,
from the   double extension method of   \cite{revoy}, \cite{medina},
any indecomposable quadratic solvable Lie algebra
is
a  \textit{double extension} of a quadratic solvable  Lie algebra  by $\Cmath.$

\end{example}

\begin{example}
\rm
Among the 170 (non isomorphic) nilpotent complex Lie algebras of dimension
$\leqslant 7,$
only a few are not
$\mathcal{I}$-null.
Those are listed in
 Table 1
in the classification of \cite{hindawi}, \cite{Magnin1}
(they are all $\mathcal{I}$-exact).
 Table 1 gives for each of them
$\dim \left(S^2 \mathfrak{g}^*\right) ^{\mathfrak{g}},$
 a basis for $ \left(\left(S^2 \mathfrak{g}^*\right) ^{\mathfrak{g}}
\left/ \right. \ker{\mathcal{I}}\right)$
(which in those cases is one-dimensional),
and the corresponding $I_Bs$.
$\blacksquare$ denotes quadratic Lie algebras;
for $\omega, \pi \in
\mathfrak{g}^*$,
$\circ $
stands for  the symmetric product
 $\omega \circ \pi =
 \omega \otimes \pi +
 \pi \otimes  \omega;$
 $\omega^{i,j,k}$ stands for  $ \omega^i\wedge   \omega^j\wedge  \omega^k.$
\end{example}
\begin{table}[h]%
\caption{}
\begin{flushleft}
{\fontsize{8}{7.5} \selectfont
\begin{tabular}{||p{2.2cm}@{\;}|@{\;}p{1.8cm}@{}|@{\;}p{5.3cm}@{}|@{\;}p{2.1cm}||}\hline
\hline
\textbf{algebra}&
\textbf{
$\dim \left(S^2 \mathfrak{g}^*\right) ^{\mathfrak{g}}$
}&
 basis for $ \left(S^2 \mathfrak{g}^*\right) ^{\mathfrak{g}}
\left/ \right. \ker{\mathcal{I}}$
&
$I_B$.
\\ \hline

$\mathfrak{g}_{5,4}\, \blacksquare$
&4
&$\omega^1 \odot \omega^5
-\omega^2 \odot \omega^4
+\omega^3 \otimes \omega^3$
&$\omega^{1,2,3}=d\omega^{1,5}$\\ \hline

$\mathfrak{g}_{6,3} \,\blacksquare$
&7
&$\omega^1 \odot \omega^6
-\omega^2 \odot \omega^5
+\omega^3 \odot \omega^4$
&$\omega^{1,2,3}=d\omega^{1,6}$\\ \hline

$\mathfrak{g}_{6,14}$
&4
&$\omega^1 \odot \omega^6
-\omega^2 \odot \omega^4
+\omega^3 \otimes \omega^3$
&$\omega^{1,2,3}=-d\omega^{1,4}$
\\ \hline

$\mathfrak{g}_{5,4} \times \Cmath$ \, $\blacksquare$
&7
&$\omega^1 \odot \omega^5
-\omega^2 \odot \omega^4
+\omega^3 \otimes \omega^3$
&$\omega^{1,2,3}=d\omega^{1,5}$\\ \hline

$\mathfrak{g}_{5,4} \times \Cmath^2$ \, $\blacksquare$
&11
&$\omega^1 \odot \omega^5
-\omega^2 \odot \omega^4
+\omega^3 \otimes \omega^3$
&$\omega^{1,2,3}=d\omega^{1,5}$\\ \hline

$\mathfrak{g}_{6,3} \times \Cmath$ \, $\blacksquare$
&11
&$\omega^1 \odot \omega^6
-\omega^2 \odot \omega^5
+\omega^3 \odot \omega^4
$
&$\omega^{1,2,3}=d\omega^{1,6}$\\ \hline

$\mathfrak{g}_{7,0.4(\lambda)},$
\linebreak[4]
$\mathfrak{g}_{7,0.5},$
$\mathfrak{g}_{7,0.6},$
$\mathfrak{g}_{7,1.02},$
$\mathfrak{g}_{7,1.10},$
$\mathfrak{g}_{7,1.13},$
$\mathfrak{g}_{7,1.14},$
$\mathfrak{g}_{7,1.17}$
&4
&$\omega^1 \odot \omega^5
-\omega^2 \odot \omega^4
+\omega^3 \otimes \omega^3$
&$\omega^{1,2,3}=d\omega^{1,5}$\\ \hline

$\mathfrak{g}_{7,1.03}$
&4
&$\omega^1 \odot \omega^6
-\omega^2 \odot \omega^4
+\omega^3 \otimes \omega^3$
&$\omega^{1,2,3}=d\omega^{1,6}$\\ \hline

$\mathfrak{g}_{7,2.2}$
&7
&$\omega^1 \odot \omega^4
-\omega^2 \odot \omega^6
+\omega^3 \odot \omega^5$
&$\omega^{1,2,3}=d\omega^{1,4}$\\ \hline

$\mathfrak{g}_{7,2.4}\, \blacksquare$
&4
&$\omega^1 \odot \omega^7
+\omega^2 \odot \omega^6
-\omega^3 \odot \omega^5
+\omega^4 \otimes \omega^4$
&$\omega^{1,3,4}-\omega^{1,2,5}=d\omega^{1,7}$\\ \hline

$\mathfrak{g}_{7,2.5}, $
$\mathfrak{g}_{7,2.6},$
$\mathfrak{g}_{7,2.7},$
$\mathfrak{g}_{7,2.8},$
$\mathfrak{g}_{7,2.9},$
&4
&$\omega^1 \odot \omega^5
-\omega^2 \odot \omega^4
+\omega^3 \otimes \omega^3$
&$\omega^{1,2,3}=d\omega^{1,5}$\\ \hline

$\mathfrak{g}_{7,2.18}$
&7
&$\omega^1 \odot \omega^6
-\omega^2 \odot \omega^5
+\omega^4 \otimes \omega^4$
&$\omega^{1,2,4}=d\omega^{1,6}$\\ \hline

$\mathfrak{g}_{7,2.44},$
$\mathfrak{g}_{7,3.6}$
&7
&$\omega^1 \odot \omega^6
-\omega^2 \odot \omega^5
+\omega^3 \odot \omega^4
$
&$\omega^{1,2,3}=d\omega^{1,6}$\\ \hline

$\mathfrak{g}_{7,3.23}$
&7
&$\omega^1 \odot \omega^6
-\omega^2 \odot \omega^5
+\omega^3 \otimes \omega^3
$
&$\omega^{1,2,3}=d\omega^{1,6}$\\ \hline

\end{tabular}
}
\end{flushleft}
\end{table}
\begin{remark}
\rm
There are nilpotent Lie algebras of higher dimension with
\linebreak[4]
   $ \dim \left(S^2 \mathfrak{g}^*\right) ^{\mathfrak{g}}
\left/ \right. \ker{\mathcal{I}} >1.$
For example, in the case of the $10$ dimensional Lie algebra
$\mathfrak{g}$ with commutation relations
$[x_1,x_2]=x_5,
[x_1,x_3]=x_6,
[x_1,x_4]=x_7,
[x_2,x_3]=x_8,
[x_2,x_4]=x_9,
[x_3,x_4]=x_{10},$
   $ \dim \left(S^2 \mathfrak{g}^*\right) ^{\mathfrak{g}}
\left/ \right. \ker{\mathcal{I}} = 4,$
and in the analogous case
of the 15 dimensional nilpotent
Lie algebra with $5$ generators one has
   $ \dim \left(S^2 \mathfrak{g}^*\right) ^{\mathfrak{g}}
\left/ \right. \ker{\mathcal{I}} = 10.$
Those algebras are  $\mathcal{I}$-exact and not quadratic.
\end{remark}
\begin{remark}
\rm
In the \textit{transversal to dimension} approach to the classification problem of nilpotent
Lie algebras initiated in \cite{santha1}, one first associates a generalized Cartan matrix
(abbr. GCM) $A$ to any nilpotent finite dimensional complex Lie algebra $\frak{g},$
and then looks at $\frak{g}$ as the quotient $\hat{\frak{g}}(A)_+/\frak{I}$
of the nilradical
of the Borel subalgebra
of the Kac-Moody Lie algebra
$\hat{\frak{g}}(A)$ associated to $A$ by some ideal $\frak{I}.$
Then one gets for any GCM $A$ the subproblem of classifying
(up to the action of a certain group)
all ideals of
$\hat{\frak{g}}(A)_+,$
thus getting all nilpotent Lie algebras of type $A$
(see
\cite{favre-santha},
\cite{fernandez},
\cite{fernandez-nunez},
\cite{santha4}, and the references therein).
Any indecomposable GCM is of exactly one of the 3 types \textit{finite}, \textit{affine},
\textit{indefinite} (among that last  the hyperbolic GCMs, with the property that
any connected proper subdiagram of the Dynkin diagram is of finite or affine type)
(\cite{kac},\cite{wanze}).
From corollary \ref{corollarynilradical},
the nilpotent Lie algebras that are not $\mathcal{I}$-null
all come from affine or indefinite types.
Unfortunately, that is the case of many nilpotent Lie algebras, see table 2.
\end{remark}

\begin{table}[h]%
\label{KMtable1}
\caption{Kac-Moody types for indecomposable nilpotent Lie algebras of dimension $\leq 7.$
Notations for indefinite hyperbolic are those of \cite{wanze}.
}
\begin{flushleft}
{\fontsize{8}{7.5} \selectfont
\begin{tabular}{||p{1.4cm}|p{2cm}|p{1.5cm}|p{1.5cm}|p{1.5cm}|p{1.5cm}|}\hline
\hline
\textbf{algebra}&
\textbf{GCM}&
Finite&
Affine&
Indefinite Hyperbolic&
Indefinite Not Hyperbolic
\\ \hline

$\mathfrak{g}_{3}$
&$\left( \begin{smallmatrix}
2&-1\\-1&2
\end{smallmatrix}\right)$
&         $A_1$             
&                      
&                     
&                      
\\

$\mathfrak{g}_{4}$
&$\left( \begin{smallmatrix}
2&-2\\-1&2
\end{smallmatrix}\right)$
&         $C_2$             
&                      
&                     
&                      
\\

$\mathfrak{g}_{5,1}$
&$\left( \begin{smallmatrix}
2&0&-1&0\\0&2&0&-1\\-1&0&2&0\\0&-1&0&2
\end{smallmatrix}\right)$
&         $A_2 \times A_2$             
&                      
&                     
&                      
\\

$\mathfrak{g}_{5,2}$
&$\left( \begin{smallmatrix}
2&-1&-1\\-1&2&0\\-1&0&2
\end{smallmatrix}\right)$
&         $A_3$             
&                      
&                     
&                      
\\

$\mathfrak{g}_{5,3}$
&$\left( \begin{smallmatrix}
2&-2&0\\-1&2&-1\\0&-1&2
\end{smallmatrix}\right)$
&         $B_3$             
&                      
&                     
&                      
\\

$\mathfrak{g}_{5,4}$
&$\left( \begin{smallmatrix}
2&-2\\-2&2
\end{smallmatrix}\right)$
&                      
&         $A_1^{(1)}$             
&                     
&                      
\\

$\mathfrak{g}_{5,5}$
&$\left( \begin{smallmatrix}
2&-3\\-1&2
\end{smallmatrix}\right)$
&         $G_2$             
&                      
&                     
&                      
\\

$\mathfrak{g}_{5,6}$
&$\left( \begin{smallmatrix}
2&-3\\-2&2
\end{smallmatrix}\right)$
&                      
&                      
&         $(3,2)$            
&                      
\\

$\mathfrak{g}_{6,1}$
&$\left( \begin{smallmatrix}
2&-1&0&-1\\-1&2&-1&0\\0&-1&2&0\\-1&0&0&2
\end{smallmatrix}\right)$
&         $D_4$             
&                      
&                     
&                      
\\

$\mathfrak{g}_{6,2}$
&$\left( \begin{smallmatrix}
2&-2&0&0\\-1&2&-1&0\\0&0&2&-1\\0&0&-1&2
\end{smallmatrix}\right)$
&         $B_2 \times A_2$             
&                      
&                     
&                      
\\

$\mathfrak{g}_{6,3}$
&$\left( \begin{smallmatrix}
2&-1&-1\\-1&2&-1\\-1&-1&2
\end{smallmatrix}\right)$
&                     
&        $A_2^{(1)}$              
&                     
&                      
\\

$\mathfrak{g}_{6,4}$
&$\left( \begin{smallmatrix}
2&-1&-1\\-2&2&0\\-1&0&2
\end{smallmatrix}\right)$
&       $B_3$              
&                      
&                     
&                      
\\

$\mathfrak{g}_{6,5}$
&$\left( \begin{smallmatrix}
2&-2&-1\\-2&2&-1\\-1&-1&2
\end{smallmatrix}\right)$
&                     
&                      
&       $H_2^{(3)}$              
&                      
\\

$\mathfrak{g}_{6,6}$
&$\left( \begin{smallmatrix}
2&-1&0\\-2&2&-1\\0&-1&2
\end{smallmatrix}\right)$
&       $C_3$              
&                      
&                     
&                      
\\

$\mathfrak{g}_{6,7}$
&$\left( \begin{smallmatrix}
2&-2&-1\\-1&2&-1\\-1&-1&2
\end{smallmatrix}\right)$
&                     
&                      
&       $H_1^{(3)}$              
&                      
\\

$\mathfrak{g}_{6,8}$
&$\left( \begin{smallmatrix}
2&-2&0\\-2&2&-1\\0&-1&2
\end{smallmatrix}\right)$
&                     
&                      
&       $H_{96}^{(3)}$              
&                      
\\

$\mathfrak{g}_{6,9}$
&$\left( \begin{smallmatrix}
2&-1&-1\\-1&2&0\\-1&0&2
\end{smallmatrix}\right)$
&       $A_3$              
&                      
&                     
&                      
\\

$\mathfrak{g}_{6,10}$
&$\left( \begin{smallmatrix}
2&-2&-1\\-1&2&0\\-2&0&2
\end{smallmatrix}\right)$
&                     
&       $A_4^{(2)}$               
&                     
&                      
\\

$\mathfrak{g}_{6,11}$
&$\left( \begin{smallmatrix}
2&-3&0\\-1&2&-1\\0&-1&2
\end{smallmatrix}\right)$
&                     
&       $G_2^{(1)}$               
&                     
&                      
\\
$\mathfrak{g}_{6,12}$
&$\left( \begin{smallmatrix}
2&-3&-2\\-2&2&-1\\-1&-1&2
\end{smallmatrix}\right)$
&                     
&                      
&      $H_{10}^{(3)}$               
&                      
\\

$\mathfrak{g}_{6,13}$
&$\left( \begin{smallmatrix}
2&-1&0\\-3&2&-1\\0&-1&2
\end{smallmatrix}\right)$
&                     
&       $G_2^{(1)}$               
&                     
&                      
\\

$\mathfrak{g}_{6,14}$
&$\left( \begin{smallmatrix}
2&-3\\-2&2
\end{smallmatrix}\right)$
&                     
&                      
&       $(3,2)$              
&                      
\\

$\mathfrak{g}_{6,15}$
&$\left( \begin{smallmatrix}
2&-2\\-2&2
\end{smallmatrix}\right)$
&                     
&       $A_1^{(1)}$               
&                     
&                      
\\

$\mathfrak{g}_{6,16}$
&$\left( \begin{smallmatrix}
2&-4\\-1&2
\end{smallmatrix}\right)$
&                     
&       $A_2^{(2)}$               
&                     
&                      
\\
$\mathfrak{g}_{6,17}$
&$\left( \begin{smallmatrix}
2&-4\\-2&2
\end{smallmatrix}\right)$
&                     
&                      
&       $(4,2)$              
&                      
\\
\\ \hline

\end{tabular}
}
\end{flushleft}
\end{table}

\setcounter{table}{1}
\begin{table}[h]%
\caption{continued}
\begin{flushleft}
{\fontsize{8}{7.5} \selectfont
\begin{tabular}{||p{1.4cm}|p{2cm}|p{1.5cm}|p{1.5cm}|p{1.5cm}|p{1.5cm}|}\hline
\hline
\textbf{algebra}&
\textbf{GCM}&
Finite&
Affine&
Indefinite Hyperbolic&
Indefinite Not Hyperbolic
\\ \hline

$\mathfrak{g}_{6,18}$
&$\left( \begin{smallmatrix}
2&-3\\-1&2
\end{smallmatrix}\right)$
&         $G_2$             
&                      
&                     
&                      
\\

$\mathfrak{g}_{6,19}$
&$\left( \begin{smallmatrix}
2&-4\\-1&2
\end{smallmatrix}\right)$
&                     
&       $A_2^{(2)}$               
&                     
&                      
\\

$\mathfrak{g}_{6,20}$
&$\left( \begin{smallmatrix}
2&-3\\-3&2
\end{smallmatrix}\right)$
&                     
&                      
&      $(3,3)$               
&                     
\\

$\mathfrak{g}_{7,0.1}$
&$\left( \begin{smallmatrix}
2&-5\\-5&2
\end{smallmatrix}\right)$
&                      
&                      
& $(5,5)$                    
&                      
\\

$\mathfrak{g}_{7,0.2}$
&idem
&                      
&                      
& idem                    
\\

$\mathfrak{g}_{7,0.3}$
&idem
&                      
&                      
& idem                    
\\

$\mathfrak{g}_{7,0.4(\lambda)}$
&$\left( \begin{smallmatrix}
2&-4\\-4&2
\end{smallmatrix}\right)$
&                      
&                      
& $(4,4)$                    
&                      
\\
$\mathfrak{g}_{7,0.5}$
&idem
&                      
&                      
&       idem              
&                      
\\
$\mathfrak{g}_{7,0.6}$
&$\left( \begin{smallmatrix}
2&-3\\-3&2
\end{smallmatrix}\right)$
&                      
&                      
& $(3,3)$                    
&                      
\\
$\mathfrak{g}_{7,0.7}$
&idem
&                      
&                      
&   idem                  
&                      
\\

$\mathfrak{g}_{7,0.8}$
&$\left( \begin{smallmatrix}
2&-3&-3\\-3&2&-3\\-3&-3&2
\end{smallmatrix}\right)$
&                      
&                      
&                     
& $\surd$                     
\\

$\mathfrak{g}_{7,1.01(i)}$
&$\left( \begin{smallmatrix}
2&0&-4\\0&2&-4\\-1&-1&2
\end{smallmatrix}\right)$
&                      
&                      
&  $H_{123}^{(3)}$                  
&                      
\\

$\mathfrak{g}_{7,1.01(ii)}$
&idem
&                      
&                      
&  idem                  
&                      
\\

$\mathfrak{g}_{7,1.02}$
&$\left( \begin{smallmatrix}
2&-2\\-3&2
\end{smallmatrix}\right)
$
&                      
&                      
&$(3,2)$                    
&                      
\\

$\mathfrak{g}_{7,1.03}$
&$\left( \begin{smallmatrix}
2&-3\\-2&2
\end{smallmatrix}\right)
$
&                      
&                      
&  $(3,2)$                  
&                      
\\

$\mathfrak{g}_{7,1.1(i_\lambda)}$
$\lambda \neq 0$
&$\left( \begin{smallmatrix}
2&-5\\-3&2
\end{smallmatrix}\right)
$
&                      
&                      
& $(5,3)$                    
&                      
\\
$\mathfrak{g}_{7,1.1(i_\lambda)}$
$\lambda = 0$
&$\left( \begin{smallmatrix}
2&-5\\-2&2
\end{smallmatrix}\right)
$
&                      
&                      
& $(5,2)$                    
&                      
\\

$\mathfrak{g}_{7,1.1(ii)}$
&$\left( \begin{smallmatrix}
2&-5\\-1&2
\end{smallmatrix}\right)
$
&                      
&                     
&       $(5,1)$               
&                      
\\

$\mathfrak{g}_{7,1.1(iii)}$
&$\left( \begin{smallmatrix}
2&-4\\-3&2
\end{smallmatrix}\right)
$
&                      
&                      
& $(4,3)$                    
&                      
\\

$\mathfrak{g}_{7,1.1(iv)}$
&$\left( \begin{smallmatrix}
2&-2\\-3&2
\end{smallmatrix}\right)
$
&                      
&                      
& $(4,2)$                    
&                      
\\

$\mathfrak{g}_{7,1.1(v)}$
&$\left( \begin{smallmatrix}
2&0&-4\\0&2&-2\\-2&-1&2
\end{smallmatrix}\right)
$
&                      
&                      
&                     
& $\surd$                     
\\

$\mathfrak{g}_{7,1.1(vi)}$
&$\left( \begin{smallmatrix}
2&-4&-1\\-3&2&0\\-1&0&2
\end{smallmatrix}\right)
$
&                      
&                      
&                     
& $\surd$                     
\\

$\mathfrak{g}_{7,1.2(i_\lambda)}$
&$\left( \begin{smallmatrix}
2&-3&-2\\-3&2&-2\\-1&-1&2
\end{smallmatrix}\right)
$
&                      
&                     
&                    
& $\surd$                     
\\

$\mathfrak{g}_{7,1.2(ii)}$
&$\left( \begin{smallmatrix}
2&-3&-2\\-3&2&-2\\-1&-1&2
\end{smallmatrix}\right)
$
&                      
&                     
&                    
& $\surd$                     
\\

$\mathfrak{g}_{7,1.2(iii)}$
&$\left( \begin{smallmatrix}
2&-3\\-2&2
\end{smallmatrix}\right)
$
&                      
&                     
&                    
& $\surd$                     
\\

$\mathfrak{g}_{7,1.2(iv)}$
&$\left( \begin{smallmatrix}
2&-3&-2\\-3&2&-2\\-1&-1&2
\end{smallmatrix}\right)
$
&                      
&                     
&                    
&  $\surd$                  
\\
\\ \hline

\end{tabular}
}
\end{flushleft}
\end{table}

\setcounter{table}{1}
\begin{table}[h]%
\caption{continued}
\begin{flushleft}
{\fontsize{8}{7.5} \selectfont
\begin{tabular}{||p{1.4cm}|p{2cm}|p{1.5cm}|p{1.5cm}|p{1.5cm}|p{1.5cm}|}\hline
\hline
\textbf{algebra}&
\textbf{GCM}&
Finite&
Affine&
Indefinite Hyperbolic&
Indefinite Not Hyperbolic
\\ \hline

$\mathfrak{g}_{7,1.3(i_\lambda)}$
&$\left( \begin{smallmatrix}
2&-3&-3\\-2&2&-1\\-2&-1&2
\end{smallmatrix}\right)
$
&                      
&                     
&                    
& $\surd$                     
\\

$\mathfrak{g}_{7,1.3(ii)}$
&idem
&                      
&                     
&                    
& idem                     
\\

$\mathfrak{g}_{7,1.3(iii)}$
&$\left( \begin{smallmatrix}
2&-3&-3\\-2&2&0\\-2&0&2
\end{smallmatrix}\right)
$
&                      
&                     
&                    
& $\surd$                     
\\

$\mathfrak{g}_{7,1.3(iv)}$
&$\left( \begin{smallmatrix}
2&-2&-2\\-2&2&-1\\-2&-1&2
\end{smallmatrix}\right)
$
&                      
&                     
&  $H_{18}^{(3)}$                  
&                      
\\
$\mathfrak{g}_{7,1.3(v)}$
&$\left( \begin{smallmatrix}
2&-3&-3&-2\\-2&2&-1&-1\\-2&-1&2&-1\\-2&-1&-1&2
\end{smallmatrix}\right)
$
&                      
&                     
&                    
& $\surd$                     
\\

$\mathfrak{g}_{7,1.4}$
&$\left( \begin{smallmatrix}
2&-5\\-2&2
\end{smallmatrix}\right)
$
&                      
&                     
& $(5,2)$                  
&                      
\\

$\mathfrak{g}_{7,1.5}$
&$\left( \begin{smallmatrix}
2&-4\\-2&2
\end{smallmatrix}\right)
$
&                      
&                     
&  $(4,2)$                 
&                      
\\

$\mathfrak{g}_{7,1.6}$
&$\left( \begin{smallmatrix}
2&-5\\-2&2
\end{smallmatrix}\right)
$
&                      
&                     
&  $(5,2)$                 
&                      
\\

$\mathfrak{g}_{7,1.7}$
&$\left( \begin{smallmatrix}
2&-2&-2\\-2&2&-1\\-1&-1&2
\end{smallmatrix}\right)
$
&                      
&                     
&  $H_{8}^{(3)}$                  
&                      
\\

$\mathfrak{g}_{7,1.8}$
&$\left( \begin{smallmatrix}
2&-3&0\\-2&2&-1\\0&-2&2
\end{smallmatrix}\right)
$
&                      
&                     
&                   
&  $\surd$                    
\\

$\mathfrak{g}_{7,1.9}$
&$\left( \begin{smallmatrix}
2&-3&-1\\-2&2&-1\\-1&-1&2
\end{smallmatrix}\right)
$
&                      
&                     
&                   
&  $\surd$                    
\\

$\mathfrak{g}_{7,1.10}$
&$\left( \begin{smallmatrix}
2&-4\\-3&2
\end{smallmatrix}\right)
$
&                      
&                     
&$(4,3)$                   
&                      
\\

$\mathfrak{g}_{7,1.11}$
&$\left( \begin{smallmatrix}
2&-4&-3\\-2&2&-1\\-1&-1&2
\end{smallmatrix}\right)
$
&                      
&                     
&                   
&  $\surd$                    
\\

$\mathfrak{g}_{7,1.12}$
&$\left( \begin{smallmatrix}
2&-4&-2\\-2&2&-1\\-1&-1&2
\end{smallmatrix}\right)
$
&                      
&                     
&                   
&  $\surd$                    
\\

$\mathfrak{g}_{7,1.13}$
&$\left( \begin{smallmatrix}
2&-4\\-2&2
\end{smallmatrix}\right)
$
&                      
&                     
& $(4,2)$                  
&                      
\\

$\mathfrak{g}_{7,1.14}$
&$\left( \begin{smallmatrix}
2&-3\\-3&2
\end{smallmatrix}\right)
$
&                      
&                     
&  $(3,3)$                 
&                     
\\

$\mathfrak{g}_{7,1.15}$
&$\left( \begin{smallmatrix}
2&-4&-3\\-2&2&-1\\-1&-1&2
\end{smallmatrix}\right)
$
&                      
&                     
&                   
&  $\surd$                    
\\

$\mathfrak{g}_{7,1.16}$
&$\left( \begin{smallmatrix}
2&-3&-2\\-2&2&-1\\-1&-1&2
\end{smallmatrix}\right)
$
&                      
&                     
&                   
&  $\surd$                    
\\

$\mathfrak{g}_{7,1.17}$
&$\left( \begin{smallmatrix}
2&-4\\-4&2
\end{smallmatrix}\right)
$
&                      
&                     
& $(4,4)$                  
&                      
\\

$\mathfrak{g}_{7,1.18}$
&$\left( \begin{smallmatrix}
2&-3&-1\\-2&2&-1\\-1&-1&2
\end{smallmatrix}\right)
$
&                      
&                     
&                   
&  $\surd$                    
\\

$\mathfrak{g}_{7,1.19}$
&$\left( \begin{smallmatrix}
2&-2&-2\\-2&2&-2\\-2&-2&2
\end{smallmatrix}\right)
$
&                      
&                     
&  $H_{71}^{(3)}$                  
&                      
\\

$\mathfrak{g}_{7,1.20}$
&$\left( \begin{smallmatrix}
2&-1&0\\-3&2&-1\\0&-1&2
\end{smallmatrix}\right)
$
&                      
& $D_4^{(3)}$                    
&                    
&                      
\\

$\mathfrak{g}_{7,1.21}$
&$\left( \begin{smallmatrix}
2&-3&-2\\-3&2&-2\\-1&-1&2
\end{smallmatrix}\right)
$
&                      
&                     
&                    
& $\surd$                     
\\

\\ \hline

\end{tabular}
}
\end{flushleft}
\end{table}

\setcounter{table}{1}
\begin{table}[h]%
\caption{continued}
\begin{flushleft}
{\fontsize{8}{7.5} \selectfont
\begin{tabular}{||p{1.4cm}|p{2cm}|p{1.5cm}|p{1.5cm}|p{1.5cm}|p{1.5cm}|}\hline
\hline
\textbf{algebra}&
\textbf{GCM}&
Finite&
Affine&
Indefinite Hyperbolic&
Indefinite Not Hyperbolic
\\ \hline

$\mathfrak{g}_{7,2.1(i_\lambda)}$
&$\left( \begin{smallmatrix}
2&-3&-1\\-1&2&-1\\-1&-1&2
\end{smallmatrix}\right)
$
&                      
&                      
& $ H^{(3)}_3$                    
&                      
\\

$\mathfrak{g}_{7,2.1(ii)}$
&$\left( \begin{smallmatrix}
2&-3&-1\\-1&2&0\\-1&0&2
\end{smallmatrix}\right)
$
&                      
& $D_4^{(3)}$                     
&                    
&                      
\\

$\mathfrak{g}_{7,2.1(iii)}$
&$\left( \begin{smallmatrix}
2&-3&0&0\\-1&2&-1&-1\\0&-1&2&0\\0&-1&0&2
\end{smallmatrix}\right)
$
&                      
&                      
&                    
& $\surd$                     
\\

$\mathfrak{g}_{7,2.1(iv)}$
&$\left( \begin{smallmatrix}
2&0&-1&-2\\0&2&-1&0\\-1&-1&2&-1\\-2&0&-1&2
\end{smallmatrix}\right)
$
&                      
&                      
&                    
&  $\surd$                    
\\

$\mathfrak{g}_{7,2.1(v)}$
&$\left( \begin{smallmatrix}
2&-2&-1\\-1&2&-1\\-1&-1&2
\end{smallmatrix}\right)
$
&                      
&                      
&  $H_1^{(3)}$                  
&                      
\\

$\mathfrak{g}_{7,2.2}$
&$\left( \begin{smallmatrix}
2&-1&-1\\-2&2&-1\\-2&-1&2
\end{smallmatrix}\right)
$
&                      
&                      
& $ H^{(3)}_7$                    
&                      
\\

$\mathfrak{g}_{7,2.3}$
&$\left( \begin{smallmatrix}
2&-5\\-2&2
\end{smallmatrix}\right)
$
&                      
&                      
&  $(5,1)$                    
&                      
\\

$\mathfrak{g}_{7,2.4}$
&$\left( \begin{smallmatrix}
2&-4\\-1&2
\end{smallmatrix}\right)
$
&                      
&  $A_2^{(2)}$                    
&                      
&                      
\\

$\mathfrak{g}_{7,2.5}$
&$\left( \begin{smallmatrix}
2&-2\\-2&2
\end{smallmatrix}\right)
$
&                      
&  $A_1^{(1)}$                    
&                      
&                      
\\

$\mathfrak{g}_{7,2.6}$
&$\left( \begin{smallmatrix}
2&-3\\-2&2
\end{smallmatrix}\right)
$
&                      
&                      
&  $(3,2)$                    
&                      
\\

$\mathfrak{g}_{7,2.7}$
&$\left( \begin{smallmatrix}
2&-4\\-2&2
\end{smallmatrix}\right)
$
&                      
&                      
&  $(4,2)$                    
&                      
\\

$\mathfrak{g}_{7,2.8}$
&$\left( \begin{smallmatrix}
2&-3\\-2&2
\end{smallmatrix}\right)
$
&                      
&                      
&  $(3,2)$                    
&                      
\\

$\mathfrak{g}_{7,2.9}$
&$\left( \begin{smallmatrix}
2&-3\\-3&2
\end{smallmatrix}\right)
$
&                      
&                      
&  $(3,3)$                    
&                      
\\

$\mathfrak{g}_{7,2.10}$
&$\left( \begin{smallmatrix}
2&-3&-1\\-1&2&0\\-1&0&2
\end{smallmatrix}\right)
$
&                      
&  $D_4^{(3)}$                    
&                      
&                      
\\

$\mathfrak{g}_{7,2.11}$
&$\left( \begin{smallmatrix}
2&-3&-2\\-2&2&-1\\-1&-1&2
\end{smallmatrix}\right)
$
&                      
&                     
&                      
&  $\surd$                    
\\

$\mathfrak{g}_{7,2.12}$
&$\left( \begin{smallmatrix}
2&-2&-2\\-2&2&0\\-2&0&2
\end{smallmatrix}\right)
$
&                      
&                      
&  $H^{(3)}_{109}$                    
&                      
\\

$\mathfrak{g}_{7,2.13}$
&$\left( \begin{smallmatrix}
2&-3&0\\-1&2&-2\\0&-1&2
\end{smallmatrix}\right)
$
&                      
&                      
&  $H^{(3)}_{100}$                    
&                      
\\

$\mathfrak{g}_{7,2.14}$
&$\left( \begin{smallmatrix}
2&-4&0\\-1&2&-1\\0&-2&2
\end{smallmatrix}\right)
$
&                      
&                      
&  $H^{(3)}_{107}$                    
&                      
\\

$\mathfrak{g}_{7,2.15}$
&$\left( \begin{smallmatrix}
2&-4&0\\-1&2&-1\\0&-1&2
\end{smallmatrix}\right)
$
&                      
&                      
&  $H^{(3)}_{97}$                    
&                      
\\

$\mathfrak{g}_{7,2.16}$
&
\text{ idem }
&                      
&                      
& \text{ idem }                    
&                      
\\

$\mathfrak{g}_{7,2.17}$
&$\left( \begin{smallmatrix}
2&-3&0\\-2&2&-1\\0&-1&2
\end{smallmatrix}\right)
$
&                      
&                      
&                      
& $\surd$                     
\\

$\mathfrak{g}_{7,2.18}$
&
\text{ idem }
&                      
&                      
&                      
&
\text{ idem }
\\

$\mathfrak{g}_{7,2.19}$
&$\left( \begin{smallmatrix}
2&-3&-1\\-2&2&0\\-1&0&2
\end{smallmatrix}\right)
$
&                      
&                      
&                      
& $\surd$                     
\\

$\mathfrak{g}_{7,2.20}$
&$\left( \begin{smallmatrix}
2&-1&-3\\-1&2&0\\-2&0&2
\end{smallmatrix}\right)
$
&                      
&                      
&                      
& $\surd$                     
\\
$\mathfrak{g}_{7,2.21}$
&$\left( \begin{smallmatrix}
2&-3&-1\\-1&2&-1\\-1&-1&2
\end{smallmatrix}\right)
$
&                      
&                      
&  $H^{(3)}_{1}$                    
&                      
\\

$\mathfrak{g}_{7,2.22}$
&$\left( \begin{smallmatrix}
2&0&-3\\0&2&-1\\-2&-1&2
\end{smallmatrix}\right)
$
&                      
&                      
&                      
& $\surd$                     
\\

$\mathfrak{g}_{7,2.23}$
&$\left( \begin{smallmatrix}
2&0&0&-2\\0&2&-2&0\\
0&-1&2&-1\\-1&0&-1&2
\end{smallmatrix}\right)
$
&                      
&  $D_4^{(2)}$                    
&                      
&                      

\\ \hline

\end{tabular}
}
\end{flushleft}
\end{table}
\setcounter{table}{1}
\begin{table}[h]%
\begin{flushleft}
\caption{ continued}
{\fontsize{8}{7.5} \selectfont
\begin{tabular}{||p{1.4cm}|p{2cm}|p{1.5cm}|p{1.5cm}|p{1.5cm}|p{1.5cm}|}\hline
\hline
\textbf{algebra}&
\textbf{GCM}&
Finite&
Affine&
Indefinite Hyperbolic&
Indefinite Not Hyperbolic
\\ \hline

$\mathfrak{g}_{7,2.24}$
&$\left( \begin{smallmatrix}
2&-3&0\\-1&2&-1\\0&-1&2
\end{smallmatrix}\right)
$
&                      
&  $G_2^{(1)}$                    
&                      
&                      
\\

$\mathfrak{g}_{7,2.25}$
&$\left( \begin{smallmatrix}
2&-3&0&-2\\-1&2&-1&0\\0&-1&2&-1\\-1&0&-1&2
\end{smallmatrix}\right)
$
&                      
&                      
&                      
& $\surd$                     
\\

$\mathfrak{g}_{7,2.26}$
&$\left( \begin{smallmatrix}
2&-1&-2\\-1&2&-1\\-2&-1&2
\end{smallmatrix}\right)
$
&                      
&                      
&  $H^{(3)}_{8}$                    
&                      
\\

$\mathfrak{g}_{7,2.27}$
&$\left( \begin{smallmatrix}
2&-2&-1&-1\\-2&2&0&-1\\-1&0&2&0\\-1&-1&0&2
\end{smallmatrix}\right)
$
&                      
&                      
&                      
& $\surd$                     
\\
$\mathfrak{g}_{7,2.28}$
&$\left( \begin{smallmatrix}
2&-1&-2&0\\-2&2&0&0\\-1&0&2&-1\\0&0&-1&2
\end{smallmatrix}\right)
$
&                      
&                      
&  $H^{(4)}_{40}$                    
&                      
\\

$\mathfrak{g}_{7,2.29}$
&$\left( \begin{smallmatrix}
2&-2&0&0\\-2&2&-1&0\\0&-1&2&-1\\0&0&-1&2
\end{smallmatrix}\right)
$
&                      
&                      
&                      
& $\surd$                     
\\

$\mathfrak{g}_{7,2.30}$
&$\left( \begin{smallmatrix}
2&-3&0&0\\-2&2&0&0\\0&0&2&-1\\0&0&-1&2
\end{smallmatrix}\right)
$
&                      
&                      
& $(3,2) \times A_1$                     
&                     
\\
$\mathfrak{g}_{7,2.31}$
&$\left( \begin{smallmatrix}
2&-3&0\\-1&2&-2\\0&-1&2
\end{smallmatrix}\right)
$
&                      
&                      
&  $H^{(3)}_{100}$                    
&                      
\\

$\mathfrak{g}_{7,2.32}$
&$\left( \begin{smallmatrix}
2&-3&-1\\-1&2&0\\-2&0&2
\end{smallmatrix}\right)
$
&                      
&                      
&  $H^{(3)}_{106}$                    
&                      
\\

$\mathfrak{g}_{7,2.33}$
&$\left( \begin{smallmatrix}
2&-3&0\\-1&2&-1\\0&-2&2
\end{smallmatrix}\right)
$
&                      
&                      
&  $H^{(3)}_{105}$                    
&                      
\\

$\mathfrak{g}_{7,2.34}$
&$\left( \begin{smallmatrix}
2&-2&-1\\-2&2&0\\-2&0&2
\end{smallmatrix}\right)
$
&                      
&                      
&  $H^{(3)}_{104}$                    
&                      
\\

$\mathfrak{g}_{7,2.35}$
&$\left( \begin{smallmatrix}
2&-1&-2\\-2&2&0\\-1&0&2
\end{smallmatrix}\right)
$
&                      
& $A_4^{(2)}$                     
&                     
&                      
\\

$\mathfrak{g}_{7,2.36}$
&$\left( \begin{smallmatrix}
2&0&-1&-1\\0&2&-1&-1\\-1&-2&2&0\\-1&-1&0&2
\end{smallmatrix}\right)
$
&                      
&                      
&                      
& $\surd$                     
\\

$\mathfrak{g}_{7,2.37}$
&$\left( \begin{smallmatrix}
2&-2&-1\\-2&2&0\\-1&0&2
\end{smallmatrix}\right)
$
&                      
&                      
&  $H^{(3)}_{96}$                    
&                      
\\

$\mathfrak{g}_{7,2.38}$
&$\left( \begin{smallmatrix}
2&-2&0&-1\\-2&2&-1&0\\0&-1&2&0\\-1&0&0&2
\end{smallmatrix}\right)
$
&                      
&                      
&                      
&  $\surd$                    
\\

$\mathfrak{g}_{7,2.39}$
&$\left( \begin{smallmatrix}
2&-2&-2\\-1&2&-1\\-1&-1&2
\end{smallmatrix}\right)
$
&                      
&                      
&  $H^{(3)}_{5}$                    
&                      
\\

$\mathfrak{g}_{7,2.40}$
&$\left( \begin{smallmatrix}
2&-2&-1\\-2&2&-1\\-1&-1&2
\end{smallmatrix}\right)
$
&                      
&                      
&  $H^{(3)}_{2}$                    
&                      
\\

$\mathfrak{g}_{7,2.41}$
&$\left( \begin{smallmatrix}
2&-2&-2\\-2&2&0\\-1&0&2
\end{smallmatrix}\right)
$
&                      
&                      
&  $H^{(3)}_{99}$                    
&                      
\\

$\mathfrak{g}_{7,2.42}$
&$\left( \begin{smallmatrix}
2&-1&-2\\-2&2&-1\\-1&-1&2
\end{smallmatrix}\right)
$
&                      
&                      
&  $H^{(3)}_{6}$                    
&                      
\\

$\mathfrak{g}_{7,2.43}$
&$\left( \begin{smallmatrix}
2&-2&-2\\-1&2&0\\-2&0&2
\end{smallmatrix}\right)
$
&                      
&                      
&  $H^{(3)}_{99}$                    
&                      
\\

$\mathfrak{g}_{7,2.44}$
&$\left( \begin{smallmatrix}
2&-1&-2\\-2&2&-1\\-1&-1&2
\end{smallmatrix}\right)
$
&                      
&                      
&  $H^{(3)}_{6}$                    
&                      
\\

$\mathfrak{g}_{7,2.45}$
&$\left( \begin{smallmatrix}
2&-2&0&-1\\-1&2&-1&-1\\0&-1&2&0\\-1&-1&0&2
\end{smallmatrix}\right)
$
&                      
&                      
&                     
&  $\surd$                    
\\
 \hline
\end{tabular}
}
\end{flushleft}
\end{table}
\setcounter{table}{1}
\begin{table}[h]%
\begin{flushleft}
\caption{continued}
{\fontsize{8}{7.5} \selectfont
\begin{tabular}{||p{1.4cm}|p{2.5cm}|p{1.5cm}|p{1.5cm}|p{1.5cm}|p{1.5cm}|}\hline
\hline
\textbf{algebra}&
\textbf{GCM}&
Finite&
Affine&
Indefinite Hyperbolic&
Indefinite Not Hyperbolic
\\ \hline

$\mathfrak{g}_{7,3.1(i_\lambda)}$
&$\left( \begin{smallmatrix}
2&-1&-1\\-1&2&-1\\-1&-1&2
\end{smallmatrix}\right) $
&                      
&  $A_2^{(1)}$                    
&                     
&                     
\\

$\mathfrak{g}_{7,3.1(iii)}$
&$\left( \begin{smallmatrix}
2&-1&-1&-1\\-1&2&0&0\\-1&0&2&0\\-1&0&0&2
\end{smallmatrix}\right) $
&  $D_4$                    
&                      
&                     
&                     
\\

$\mathfrak{g}_{7,3.2}$
&$\left( \begin{smallmatrix}
2&-3&-1\\-1&2&0\\-1&0&2
\end{smallmatrix}\right) $
&                      
&  $D_4^{(3)}$                    
&                     
&                     
\\

$\mathfrak{g}_{7,3.3}$
&$\left( \begin{smallmatrix}
2&-3&0\\-1&2&-1\\0&-1&2
\end{smallmatrix}\right) $
&                      
&  $G_2^{(1)}$                    
&                     
&                     
\\

$\mathfrak{g}_{7,3.4}$
&$\left( \begin{smallmatrix}
2&-1&-1\\-2&2&0\\-2&0&2
\end{smallmatrix}\right) $
&                      
&  $D_3^{(2)}$                    
&                     
&                     
\\

$\mathfrak{g}_{7,3.5}$
&$\left( \begin{smallmatrix}
2&-1&-1\\-2&2&0\\-1&0&2
\end{smallmatrix}\right) $
&  $B_3$                    
&                      
&                     
&                     
\\

$\mathfrak{g}_{7,3.6}$
&$\left( \begin{smallmatrix}
2&-1&-2\\-1&2&-1\\-1&-1&2
\end{smallmatrix}\right) $
&                      
&                      
&  $H^{(3)}_1$                   
&                     
\\

$\mathfrak{g}_{7,3.7}$
&$\left( \begin{smallmatrix}
2&-2&0&0\\-1&2&0&-1\\0&0&2&-1\\0&-1&-1&2
\end{smallmatrix}\right) $
&  $B_4$                    
&                      
&                     
&                     
\\

$\mathfrak{g}_{7,3.8}$
&$\left( \begin{smallmatrix}
2&-2&-1&0\\-1&2&0&-1\\-1&0&2&0\\0&-1&0&2
\end{smallmatrix}\right) $
&  $F_4$                    
&                      
&                     
&                     
\\

$\mathfrak{g}_{7,3.9}$
&$\left( \begin{smallmatrix}
2&-2&0&0\\-1&2&-1&-1\\0&-1&2&0\\0&-1&0&2
\end{smallmatrix}\right) $
&                      
&  $B_4^{(1)}$                    
&                     
&                     
\\

$\mathfrak{g}_{7,3.10}$
&$\left( \begin{smallmatrix}
2&-1&-1&0\\-1&2&0&-2\\-1&0&2&0\\0&-1&0&2
\end{smallmatrix}\right) $
&  $C_4$                    
&                      
&                     
&                     
\\

$\mathfrak{g}_{7,3.11}$
&$\left( \begin{smallmatrix}
2&-2&-1&0\\-1&2&0&-1\\-1&0&2&0\\0&-1&0&2
\end{smallmatrix}\right) $
&  $F_4$                    
&                      
&                     
&                     
\\

$\mathfrak{g}_{7,3.12}$
&$\left( \begin{smallmatrix}
2&-1&-1&0\\-1&2&0&-1\\-1&0&2&-1\\0&-1&-1&2
\end{smallmatrix}\right) $
&                      
&  $A^{(1)}_3$                    
&                     
&                     
\\

$\mathfrak{g}_{7,3.13}$
&$\left( \begin{smallmatrix}
2&-2&0&0\\-2&2&0&0\\0&0&2&-1\\0&0&-1&2
\end{smallmatrix}\right) $
&                      
&  $A^{(1)}_1 \times A_1$                    
&                     
&                     
\\

$\mathfrak{g}_{7,3.14}$
&$\left( \begin{smallmatrix}
2&-1&-2&0\\-1&2&0&-1\\-1&0&2&0\\0&-1&0&2
\end{smallmatrix}\right) $
&  $C_4$                    
&                      
&                     
&                     
\\

$\mathfrak{g}_{7,3.15}$
&$\left( \begin{smallmatrix}
2&-1&-1&0\\-2&2&0&0\\-1&0&2&-1\\0&0&-1&2
\end{smallmatrix}\right) $
&  $B_4$                    
&                      
&                     
&                     
\\

$\mathfrak{g}_{7,3.16}$
&$\left( \begin{smallmatrix}
2&-2&0&0\\-1&2&0&0\\0&0&2&-2\\0&0&-1&2
\end{smallmatrix}\right) $
&  $B_2 \times B_2$                    
&                      
&                     
&                     
\\

$\mathfrak{g}_{7,3.17}$
&$\left( \begin{smallmatrix}
2&-3&0&0\\-1&2&0&0\\0&0&2&-1\\0&0&-1&2
\end{smallmatrix}\right) $
&  $G_2 \times A_2$                    
&                      
&                     
&                     
\\

$\mathfrak{g}_{7,3.18}$
&$\left( \begin{smallmatrix}
2&-2&0&0&0\\-1&2&0&0&-1\\0&0&2&-1&0\\0&0&-1&2&0\\0&-1&0&0&2
\end{smallmatrix}\right) $
&  $B_3 \times A_2$                    
&                      
&                     
&                     
\\

$\mathfrak{g}_{7,3.19}$
&$\left( \begin{smallmatrix}
2&-1&-1&0&0\\-1&2&0&0&0\\-1&0&2&-1&0\\0&0&-1&2&-1\\0&0&0&-1&2
\end{smallmatrix}\right) $
&  $A_5$                    
&                      
&                     
&                     
\\
 \hline
\end{tabular}
}
\end{flushleft}
\end{table}
\setcounter{table}{1}
\begin{table}[h]%
\begin{flushleft}
\caption{continued}
{\fontsize{8}{7.5} \selectfont
\begin{tabular}{||p{1.4cm}|p{2.8cm}|p{1.8cm}|p{1.5cm}|p{1.5cm}|p{1.5cm}|}\hline
\hline
\textbf{algebra}&
\textbf{GCM}&
Finite&
Affine&
Indefinite Hyperbolic&
Indefinite Not Hyperbolic
\\ \hline

$\mathfrak{g}_{7,3.20}$
&$\left( \begin{smallmatrix}
2&-2&-2\\-1&2&0\\-1&0&2
\end{smallmatrix}\right) $
&                      
&  $C_2^{(1)}$                    
&                     
&                     
\\

$\mathfrak{g}_{7,3.21}$
&$\left( \begin{smallmatrix}
2&-2&-1\\-1&2&0\\-2&0&2
\end{smallmatrix}\right) $
&                      
&  $A_4^{(2)}$                    
&                     
&                     
\\

$\mathfrak{g}_{7,3.22}$
&$\left( \begin{smallmatrix}
2&-1&-2\\-1&2&0\\-1&0&2
\end{smallmatrix}\right) $
&  $C_3$                    
&                      
&                     
&                     
\\

$\mathfrak{g}_{7,3.23}$
&$\left( \begin{smallmatrix}
2&-2&-1\\-2&2&0\\-1&0&2
\end{smallmatrix}\right) $
&                      
&                      
&   $H^{(3)}_{96}$                  
&                     
\\

$\mathfrak{g}_{7,3.24}$
&$\left( \begin{smallmatrix}
2&-1&0&0\\-1&2&-1&-1\\0&-1&2&-1\\0&-1&-1&2
\end{smallmatrix}\right) $
&                      
&                      
&   $H^{(4)}_{3}$                  
&                     
\\

$\mathfrak{g}_{7,4.1}$
&$\left( \begin{smallmatrix}
2&-1&-1&0\\-1&2&0&0\\-1&0&2&-1\\0&0&-1&2
\end{smallmatrix}\right) $
&   $A_4$                   
&                      
&                     
&                     
\\

$\mathfrak{g}_{7,4.2}$
&$\left( \begin{smallmatrix}
2&-1&-1&-1\\-1&2&0&0\\-1&0&2&0\\-1&0&0&2
\end{smallmatrix}\right) $
&   $D_4$                   
&                      
&                     
&                     
\\

$\mathfrak{g}_{7,4.3}$
&$\left( \begin{smallmatrix}
2&-1&0&0&0\\-1&2&0&0&0\\0&0&2&0&-1\\0&0&0&2&-1\\0&0&-1&-1&2
\end{smallmatrix}\right) $
&   $A_4 \times A_3$                   
&                      
&                     
&                     
\\

$\mathfrak{g}_{7,4.4}$
&$\left( \begin{smallmatrix}
2&0&0&-1&0&0\\
0&2&0&0&-1&0\\
0&0&2&0&0&-1\\
-1&0&0&2&0&0\\
0&-1&0&0&2&0\\
0&0&-1&0&0&2
\end{smallmatrix}\right) $
&   $A_2 \times A_2 \times A_2$                   
&                      
&                     
&                     
\\
 \hline
\end{tabular}
}
\end{flushleft}
\end{table}
\begin{example}
\rm
The quadratic $5$-dimensional nilpotent Lie algebra
$\mathfrak{g}_{5,4}$
has commutation relations
$[x_1,x_2]=x_3,$
$[x_1,x_3]=x_4,$
$[x_2,x_3]=x_5.$
Consider the $10$-dimensional direct product
$\mathfrak{g}_{5,4} \times \mathfrak{g}_{5,4},$
with the commutation relations:
$[x_1,x_2]=x_5,$
$[x_1,x_5]=x_6,$
$[x_2,x_5]=x_7,$
$[x_3,x_4]=x_8,$
$[x_3,x_8]=x_9,$
$[x_4,x_8]=x_{10}.$
The only $11$-dimensional nilpotent Lie algebra
with an invariant bilinear form which reduces to
$B_1=\omega^1 \odot \omega^7
-\omega^2 \odot \omega^6
+\omega^5 \otimes \omega^5,$
$B_2=\omega^3 \odot \omega^{10}
-\omega^4 \odot \omega^9
+\omega^8 \otimes \omega^8,$
on respectively the first and second factor is
the direct product
$\Cmath \times \mathfrak{g}_{5,4} \times \mathfrak{g}_{5,4},$
\end{example}

\begin{example}
\rm
The 4-dimensional solvable "diamond" Lie algebra $\frak{g}$ with basis
$(x_1,x_2,x_3,x_4)$ and commutation relations
$[x_1,x_2]=x_3,
[x_1,x_3]=-x_2,
[x_2,x_3]=x_4$
cannot be obtained as in lemma
\ref{lemmeborel}.
Here
 $ \dim \left(\left(S^2 \mathfrak{g}^*\right) ^{\mathfrak{g}}
\left/ \right. \ker{\mathcal{I}}\right)=1,$
with basis element
$B= \omega^1 \circ \omega^4  + \omega^2 \otimes \omega^2
+ \omega^3 \otimes \omega^3.$
 $I_B= \omega^{1,2,3} =d\omega^{1,4};$
$\frak{g}$ is quadratic and $\mathcal{I}$-exact.
In fact, one verifies that all other solvable 4-dimensional Lie solvable Lie algebras are
 $\mathcal{I}$-null (for a list, see e.g. \cite{ovando}).

\end{example}
\par

\section{Comparison of $HL^2$ and $H^2$ for a Lie algebra.}

If $\mathfrak{c}$ denotes the center of
$\mathfrak{g},$
$\mathfrak{c} \,  \otimes \left(S^2 \mathfrak{g}^*\right) ^{\mathfrak{g}}$
is the space of
invariant $\mathfrak{c}$-valued symmetric bilinear map and
we denote
$F= Id\, \otimes \mathcal{I}
: \,
\mathfrak{c} \,  \otimes \left( S^2 \mathfrak{g}^*\right) ^{\mathfrak{g}}
\rightarrow C^3(\mathfrak{g}, \mathfrak{g})=
\mathfrak{g} \,  \otimes \bigwedge^3 \mathfrak{g}^*.$
Then
$\text{Im\,}F= \mathfrak{c}\, \otimes \text{Im\,}{\mathcal{I}}.$

\begin{theorem}
Let $\mathfrak{g}$ be any finite dimensional complex Lie algebra and
$ZL^2_0(\mathfrak{g},\mathfrak{g})$
(resp. $ZL^2_0(\mathfrak{g},\Cmath))$
the space of
symmetric adjoint (resp. trivial) Leibniz 2-cocycles.
\\(i)
$ZL^2(\mathfrak{g},\mathfrak{g}) \left/
\left(Z^2(\mathfrak{g},\mathfrak{g}) \oplus ZL_0^2(\mathfrak{g},\mathfrak{g})\right)
\right. \cong
\left(\mathfrak{c}\, \otimes \text{Im\,}{\mathcal{I}}\right) \cap B^3(\mathfrak{g},\mathfrak{g}) .$
\\(ii)
$ZL^2_0(\mathfrak{g},\mathfrak{g})= \mathfrak{c}\, \otimes \ker{\mathcal{I}}.$
In particular,
$\dim ZL^2_0(\mathfrak{g},\mathfrak{g}) = c \, \frac{p(p+1)}{2} $
where
$c=\dim  \mathfrak{c} $ and
$p=
\dim  \mathfrak{g}/\mathcal{C}^2 \mathfrak{g}
=\dim H^1( \mathfrak{g},\Cmath).$
 \smallskip
\\(iii)
$HL^2(\mathfrak{g},\mathfrak{g}) \cong H^2(\mathfrak{g},\mathfrak{g}) \oplus
\left(\mathfrak{c} \otimes \ker{\mathcal{I}}\right)
\oplus \left(
\left(\mathfrak{c} \otimes \text{Im\,}{\mathcal{I}}\right) \cap B^3(\mathfrak{g},\mathfrak{g}) \right).$
\\(iv)
$ZL^2(\mathfrak{g},\Cmath) \left/
\left(Z^2(\mathfrak{g},\Cmath) \oplus ZL_0^2(\mathfrak{g},\Cmath)\right)
\right. \cong
\text{Im\,}{\mathcal{I}} \cap B^3(\mathfrak{g},\Cmath) .$
\\(v)
$ZL^2_0(\mathfrak{g},\Cmath)= \ker{\mathcal{I}}.$
\\(vi)
$HL^2(\mathfrak{g},\Cmath) \cong H^2(\mathfrak{g},\Cmath) \oplus
\ker{\mathcal{I}}
\oplus
\left(\text{Im\,}{\mathcal{I}} \cap B^3(\mathfrak{g},\Cmath) \right).$
 \smallskip
\end{theorem}
\begin{proof}
(i) The Leibniz 2-cochain space
$CL^2(\mathfrak{g},\mathfrak{g}) =
\mathfrak{g}
\otimes
\left(\mathfrak{g}^*\right)^{\otimes 2}
$
decomposes as
$\left(
\mathfrak{g}
\otimes
\bigwedge^2 \mathfrak{g}^*  \right) \;
 \oplus \; \left(
\mathfrak{g}
\otimes
 S^2 \,\mathfrak{g}^*  \right)$
 with
 $
\mathfrak{g}
\otimes
 S^2 \,\mathfrak{g}^*  $
 the space of symmetric elements in
$CL^2(\mathfrak{g},\mathfrak{g}).$
By definition of the Leibniz coboundary $\delta$, one has for
$\psi \in CL^2(\mathfrak{g},\mathfrak{g})$ and $X,Y,Z \in \mathfrak{g}$
\begin{equation}
\label{leibnizcocycle}
(\delta \psi)(X,Y,Z) = u+v+w+r+s+t
\end{equation}
with $u=[X,\psi(Y,Z)],\;
v=[\psi(X,Z),Y], \;
w=-[\psi(X,Y),Z],\;
r=-\psi([X,Y],Z), \;
$
$s=\psi(X,[Y,Z]),  \;
t=\psi([X,Z],Y).$
$\delta$ coincides with the usal coboundary operator on
$
\mathfrak{g}
\otimes
\bigwedge^2 \mathfrak{g}^*   .$
Now, let $\psi = \psi_1 + \psi_0 \in
CL^2(\mathfrak{g},\mathfrak{g})$
, $\psi_1 \in
\mathfrak{g}
\otimes
\bigwedge^2 \mathfrak{g}^*   , \;
 \psi_0 \in
\mathfrak{g}
\otimes
 S^2 \,\mathfrak{g}^*  .$
\\
\par
 Suppose
$\psi \in ZL^2(\mathfrak{g},\mathfrak{g}):$
$\delta \psi = 0 = \delta\psi_1 + \delta \psi_0 =  d \psi_1 + \delta \psi_0.$
Then $\delta \psi_0 = -d \psi_1
\in
\mathfrak{g}
\otimes
\bigwedge^3 \mathfrak{g}^*  $ is antisymmetric.
Then permuting $X$ and $Y$ in formula
(\ref{leibnizcocycle}) for $\psi_0$ yields
$(\delta \psi_0)(Y,X,Z) = -v - u + w - r + t + s.$
As
$\delta \psi_0$ is antisymmetric, we get
\begin{equation}
\label{first}
w + s + t = 0.
\end{equation}
Now, the circular permutation $(X,Y,Z)$ in
(\ref{leibnizcocycle}) for $\psi_0$ yields
$(\delta \psi_0)(Y,Z,X) = -v - w + u - s - t + r.$
Again, by antisymmetry, we get
\begin{equation}
\label{second}
v +w + s + t = 0,
\end{equation}
i.e.
$(\delta \psi_0)(X,Y,Z) = u + r.$
From (\ref{first}) and (\ref{second}), $v = 0.$
Applying twice
the circular permutation $(X,Y,Z)$ to $v$, we get
first $w=0$ and then $u=0.$
Hence $(\delta \psi_0)(X,Y,Z) = r
=-\psi_0([X,Y],Z).$
Note first that $u=0$ reads
$[X,\psi_0(Y,Z)] = 0.$
As $X,Y,Z$ are arbitrary, $\psi_0$ is
$\mathfrak{c}$-valued.
Now the permutation  of $Y$ and $Z$ changes
$r$ to $- t=s$ (from (\ref{second})).
Again, by antisymmetry
of $\delta \psi_0,$
$r=t=-s.$
As $X,Y,Z$ are arbitrary, one gets
$\psi_0 \in \mathfrak{c} \,  \otimes \left(S^2 \mathfrak{g}^*\right) ^{\mathfrak{g}}.$
Now $F(\psi_0)= - r = -\delta\psi_0 = d \psi_1 \in
B^3(\mathfrak{g},\mathfrak{g} ).$
Hence
$$\psi_0 \in ZL^2_0(\mathfrak{g},\mathfrak{g})
\Leftrightarrow
F(\psi_0)=0
\Leftrightarrow
\psi_1 \in Z^2(\mathfrak{g},\mathfrak{g})
\Leftrightarrow
\psi_0 \in \mathfrak{c} \otimes \ker{\mathcal{I}}.$$
Consider now the linear map
$\Phi \; : \,
ZL^2(\mathfrak{g},\mathfrak{g}) \rightarrow
F^{-1}(B^3(\mathfrak{g},\mathfrak{g})) \left/\; \ker{F} \right.$
defined by
$\psi \mapsto [\psi_0] \,(\text{mod} \ker{F}).$
$\Phi$ is onto: for any $[\varphi_0] \in
F^{-1}(B^3(\mathfrak{g},\mathfrak{g})) \left/\; \ker{F} \right.$
,
$\varphi_0 \in
\mathfrak{c}\, \otimes  \left(S^2\mathfrak{g^*}\right)^{\mathfrak{g}},$
one has $F(\varphi_0) \in B^3(\mathfrak{g},\mathfrak{g}),$
hence $F(\varphi_0)= d\varphi_1,$ $\varphi_1 \in C^2(\mathfrak{g},\mathfrak{g}),$
and then $\varphi=\varphi_0 +\varphi_1$ is a Leibniz cocycle such that
$\Phi( \varphi) =[\varphi_0].$
Now $\ker{\Phi}=
Z^2(\mathfrak{g},\mathfrak{g}) \oplus ZL_0^2(\mathfrak{g},\mathfrak{g}),$
since condition $[\psi_0] =[0]$ reads $\psi_0 \in \ker{F}$ which is equivalent
to $\psi \in
Z^2(\mathfrak{g},\mathfrak{g}) \oplus ZL_0^2(\mathfrak{g},\mathfrak{g}).$
Hence $\Phi$ yields an isomorphism
$ZL^2(\mathfrak{g},\mathfrak{g}) \left/
\left(Z^2(\mathfrak{g},\mathfrak{g}) \oplus ZL_0^2(\mathfrak{g},\mathfrak{g})\right)
\right. \cong
F^{-1}(B^3(\mathfrak{g},\mathfrak{g})) \left/ \;\ker{F} \right. .$
The latter is isomorphic to
$ \text{Im\,}F  \cap B^3(\mathfrak{g},\mathfrak{g})
\cong
\left(\mathfrak{c} \otimes \text{Im\,}{\mathcal{I}}\right) \cap B^3(\mathfrak{g},\mathfrak{g}).$
\\(ii)  Results from the invariance of
 $\psi_0 \in ZL^2_0(\mathfrak{g},\mathfrak{g}).$
\\(iii)
Results immediately from (i), (ii) since
 $BL^2(\mathfrak{g},\mathfrak{g})=B^2(\mathfrak{g},\mathfrak{g})$
 as the Leibniz differential on
$CL^1(\mathfrak{g},\mathfrak{g}) = \mathfrak{g}^* \otimes \mathfrak{g} = C^1(\mathfrak{g},\mathfrak{g})$
coincides with the usual one.
\\(iv)-(vi)
Readily similar.
\end{proof}

\begin{remark}
\rm
Since
$ \ker{\mathcal{I}}\oplus \left(\text{Im\,}{\mathcal{I}} \cap B^3(\mathfrak{g},\Cmath) \right) \cong \ker{h}$ where $h$ denotes $\mathcal{I}$ composed with the projection of
$Z^3(\mathfrak{g},\Cmath)$
onto $H^3(\mathfrak{g},\Cmath),$ the result (vi) is the same as
in \cite{hpl}.
\end{remark}

\begin{remark}
\rm
Any supplementary subspace to
$Z^2(\mathfrak{g},\Cmath) \oplus ZL_0^2(\mathfrak{g},\Cmath)$ in
$ZL^2(\mathfrak{g},\Cmath)$
consists of \textit{coupled} Leibniz 2-cocycles, i.e. the nonzero elements
have the property that their
symmetric and antisymmetric
parts are not cocycles.
To get such a supplementary subspace, pick any supplementary subspace $W$ to
$\ker{\mathcal{I}}$ in
$\left(S^2 \mathfrak{g}^*\right) ^{\mathfrak{g}}$
and take $\mathcal{C} =\left\{ B +\omega\, ; B \in W \cap {\mathcal{I}}^{-1}
(B^3(\mathfrak{g},\Cmath)), I_B=d\omega \right\}.$
\end{remark}
\begin{definition}
$\mathfrak{g}$ is said to be adjoint (resp. trivial) $ZL^2$-uncoupling
if
\linebreak[4]

$\left(\mathfrak{c}\, \otimes \text{Im\,}{\mathcal{I}}\right) \cap B^3(\mathfrak{g},\mathfrak{g}) =\{0\}$
(resp. $\text{Im\,}{\mathcal{I}} \cap B^3(\mathfrak{g},\Cmath) =\{0\}).$
\end{definition}
Adjoint $ZL^2$-uncoupling implies trivial $ZL^2$-uncoupling,
since
\linebreak[4]
$\mathfrak{c}\, \otimes \left( \text{Im\,}{\mathcal{I}} \cap B^3(\mathfrak{g},\Cmath)\right) \subset
\left(\mathfrak{c}\, \otimes \text{Im\,}{\mathcal{I}}\right) \cap B^3(\mathfrak{g},\mathfrak{g}).$
The reciprocal holds obviously true
for $\mathcal{I}$-exact or zero-center Lie algebras.
However we do not know
if it holds true in general
(e.g. we do not know of a nilpotent Lie algebra which is not
 $\mathcal{I}$-exact).
The class of adjoint
$ZL^2$-uncoupling  Lie algebras is rather extensive since it contains all zero-center Lie algebras
and all
$\mathcal{I}$-null Lie algebras.
\begin{corollary}
(i)
$HL^2(\mathfrak{g},\mathfrak{g}) \cong H^2(\mathfrak{g},\mathfrak{g}) \oplus
\left(\mathfrak{c} \otimes \ker{\mathcal{I}}\right)$
if and only if
$\mathfrak{g}$ is adjoint $ZL^2$-uncoupling.
\\(ii)
$HL^2(\mathfrak{g},\Cmath) \cong H^2(\mathfrak{g},\Cmath) \oplus
\ker{\mathcal{I}}$
if and only if
$\mathfrak{g}$ is trivial $ZL^2$-uncoupling.
\end{corollary}
\par
\begin{example}
\rm
$\mathfrak{g} ={\mathfrak{g}}_{5,4}.$
\\
Trivial Leibniz cohomology.
$B^2(\mathfrak{g},\Cmath) = \langle
d\omega^3=-\omega^{1,2}, d\omega^4=-\omega^{1,3}, d\omega^5=-\omega^{2,3}\rangle,$
$\dim Z^2(\mathfrak{g},\Cmath) = 6,$
$\dim H^2(\mathfrak{g},\Cmath) = 3,$
$Z^2(\mathfrak{g},\Cmath) = \langle
\omega^{1,4},
\omega^{2,5},
\omega^{1,5}+
\omega^{2,4}\rangle \oplus  B^2(\mathfrak{g},\Cmath),$
$\dim{ZL^2_0(\mathfrak{g},\Cmath)}=3,$
$ZL^2_0(\mathfrak{g},\Cmath)
(\cong \ker{\mathcal{I}})
= \langle
\omega^1 \otimes \omega^1,
\omega^1 \odot \omega^2,
\omega^2 \otimes \omega^2 \rangle,$
$\dim{ZL^2(\mathfrak{g},\Cmath) }= 10,$
$\dim{HL^2(\mathfrak{g},\Cmath) }= 7,$ and
\begin{eqnarray*}
ZL^2(\mathfrak{g},\Cmath)& =&
Z^2(\mathfrak{g},\Cmath) \oplus
ZL^2_0(\mathfrak{g},\Cmath)
\oplus \Cmath g_1,\\
HL^2(\mathfrak{g},\Cmath) &=&
H^2(\mathfrak{g},\Cmath) \oplus
ZL^2_0(\mathfrak{g},\Cmath)
\oplus \Cmath g_1
\end{eqnarray*}
with
$g_1= B + \omega^{1,5}$
and $B= \omega^1 \odot \omega^5 -  \omega^2 \odot \omega^4  +\omega^3 \otimes \omega^3,$
(we already know that
$\text{Im\,}{\mathcal{I}} = \Cmath I_B = \Cmath d\omega^{1,5}$
and
$\text{Im\,}{\mathcal{I}} \cap B^3(\mathfrak{g},\Cmath)=
\text{Im\,}{\mathcal{I}} $
is one-dimensional).
${\mathfrak{g}}_{5,4}$  is not trivial $ZL^2$-uncoupling
(hence not adjoint $ZL^2$-uncoupling either),
and $g_1$ is a coupled Leibniz 2-cocycle.
\\Adjoint Leibniz cohomology.
$\dim{Z^2(\mathfrak{g},\mathfrak{g})} = 24;$
$ZL^2_0(\mathfrak{g},\mathfrak{g})= \mathfrak{c}\, \otimes \ker{\mathcal{I}}$ has dimension 6,
$\dim{ZL^2(\mathfrak{g},\mathfrak{g})}=32,$
\begin{eqnarray*}
ZL^2(\mathfrak{g},\mathfrak{g})&=&
Z^2(\mathfrak{g},\mathfrak{g})
\oplus
ZL^2_0(\mathfrak{g},\mathfrak{g})
\oplus \Cmath G_1
\oplus \Cmath G_2,\\
HL^2(\mathfrak{g},\mathfrak{g})&=&
H^2(\mathfrak{g},\mathfrak{g})
\oplus
ZL^2_0(\mathfrak{g},\mathfrak{g})
\oplus \Cmath G_1
\oplus \Cmath G_2,
\end{eqnarray*}
where $G_1,G_2$ are the following Leibniz 2-cocycles, each of which is coupled:
\begin{eqnarray*}
G_1&=&
x_5 \otimes ( B +  \omega^{1,5})\\
G_2&=&
x_4 \otimes ( B +   \omega^{1,5})
.
\end{eqnarray*}

\end{example}

\par
\begin{corollary}
For any Lie algebra
$\mathfrak{g}$ with trivial center
$\mathfrak{c}= \{0\},$
 $HL^2(\mathfrak{g},\mathfrak{g})=H^2(\mathfrak{g},\mathfrak{g}).$
 In particular any rigid
Lie algebra  with trivial center is Leibniz rigid.
\end{corollary}

\begin{corollary}
For any
reductive algebra
Lie
$\mathfrak{g}$ with  center
$\mathfrak{c},$
$HL^2(\mathfrak{g},\mathfrak{g}) =
H^2(\mathfrak{g},\mathfrak{g}) \oplus \,
 \left(\mathfrak{c} \, \otimes S^2 \mathfrak{c}^* \right),$
and
$\dim H^2(\mathfrak{g},\mathfrak{g}) =  \frac{c^2(c-1)}{2}$
 with $c=\dim \mathfrak{c}.$
\end{corollary}
\begin{proof}
$\mathfrak{g} = \mathfrak{s} \oplus \mathfrak{c}$
with $\mathfrak{s} = \mathcal{C}^2  \mathfrak{g}$ semisimple.
We first prove that
$\mathfrak{g}$ is adjoint $ZL^2$-uncoupling.
$ \mathfrak{c} \,  \otimes \left(S^2 \mathfrak{g}^*\right)^{\mathfrak{g}}
= \mathfrak{c} \,  \otimes \left(S^2 \mathfrak{s}^*\right) ^{\mathfrak{s}}
\oplus \left(\mathfrak{c} \,  \otimes S^2 \mathfrak{c}^*\right)
= c \, \left( S^2 \mathfrak{s}^*\right) ^{\mathfrak{s}}
\oplus c \, \left( S^2 \mathfrak{c}^*\right).$
Suppose first
 $\mathfrak{s} $
 simple. Then any bilinear symmetric invariant form on
 $\mathfrak{s} $   is some multiple of the Killing form $K.$
 Hence
$ \mathfrak{c} \,  \otimes \left(S^2 \mathfrak{g}^*\right)^{\mathfrak{g}}
= c \, (\Cmath K)
\oplus c \, \left( S^2 \mathfrak{c}^*\right).$
For any $\psi_0 \in  \mathfrak{c} \,  \otimes\left( S^2 \mathfrak{g}^*\right)^{\mathfrak{g}},$
$F({\psi_0})$ is then some linear combination of copies of $I_K.$
As is well-known, $I_K$  is no coboundary. Hence  if we suppose that
$F({\psi_0})$ is a coboundary, necessarily
$F({\psi_0})=0.$
$\mathfrak{g}$ is adjoint $ZL^2$-uncoupling
when $\mathfrak{s} $  is  simple.
Now, if
 $\mathfrak{s} $  is not simple,
 $\mathfrak{s} $ can be decomposed as a direct sum
 $ \mathfrak{s}_1 \oplus \cdots \oplus   \mathfrak{s}_m$ of simple ideals of
 $\mathfrak{s}. $
 Then
$  \left( S^2 \mathfrak{s}^*\right) ^{\mathfrak{s}} =
\bigoplus_{i=1}^{m} \,    \left( S^2 \mathfrak{s_i}^*\right) ^{\mathfrak{s_i}}
=
\bigoplus_{i=1}^{m} \,    \Cmath \, K_i$
($K_i $ Killing form of
 $\mathfrak{s}_i. $)
The same reasoning then applies and shows that
$\mathfrak{g}$ is adjoint $ZL^2$-uncoupling.
From (ii) in the Theorem,
$ZL^2_0(\mathfrak{g},\mathfrak{g}) =  \mathfrak{c} \, \otimes S^2 \mathfrak{c}^* .$
Now,
$\mathfrak{g} = \mathfrak{s} \oplus \mathfrak{c}$
with $\mathfrak{s} = \mathcal{C}^2  \mathfrak{g}$ semisimple.
 $\mathfrak{s} $ can be decomposed as a direct sum
 $ \mathfrak{s}_1 \oplus \cdots \oplus   \mathfrak{s}_m$ of ideals of
 $\mathfrak{s} $
 hence of $\mathfrak{g}. $  Then
$H^2(\mathfrak{g},\mathfrak{g}) =
\bigoplus_{i=1}^{m} H^2(\mathfrak{g},\mathfrak{s}_i)\,  \oplus \,  H^2(\mathfrak{g},\mathfrak{c}).$
As
 $\mathfrak{s}_i $ is a nontrivial
$\mathfrak{g}$-module,
$H^2(\mathfrak{g},\mathfrak{s}_i) = \{0\}$ (\cite{guichardet}, Prop. 11.4, page 154).
Hence
$H^2(\mathfrak{g},\mathfrak{g}) =
H^2(\mathfrak{g},\mathfrak{c}) = c\, H^2(\mathfrak{g},\Cmath).$
By the K\"{u}nneth formula and
Whitehead's lemmas,
$H^2(\mathfrak{g},\Cmath) =
\left(H^2(\mathfrak{s},\Cmath) \otimes H^0(\mathfrak{c},\Cmath) \right) \, \oplus
\left(H^1(\mathfrak{s},\Cmath) \otimes H^1(\mathfrak{c},\Cmath) \right) \, \oplus
\left(H^0(\mathfrak{s},\Cmath) \otimes H^2(\mathfrak{c},\Cmath) \right)
= $
$H^0(\mathfrak{s},\Cmath) \otimes H^2(\mathfrak{c},\Cmath)
= \Cmath \otimes H^2(\mathfrak{c},\Cmath). $
Hence $\dim H^2(\mathfrak{g},\mathfrak{g}) =  \frac{c^2(c-1)}{2}.$
\end{proof}
\begin{example}
\rm
For
$\mathfrak{g} = \mathfrak{gl}(n),$
$HL^2(\mathfrak{g},\mathfrak{g}) =
ZL^2_0(\mathfrak{g},\mathfrak{g}) = \Cmath \,\left( x_{n^2} \oplus (
\omega^{n^2} \circ \omega^{n^2} )\right),$  where
$(x_i)_{1\leqslant i \leqslant n^2}$ is a basis of
$\mathfrak{g}$
such that
$(x_i)_{1\leqslant i \leqslant n^2-1}$ is a basis of
$ \mathfrak{sl}(n)$ and $x_{n^2}$ is the identity matrix, and
$(\omega_i)_{1\leqslant i \leqslant n^2}$ the dual basis
to $(x_i)_{1\leqslant i \leqslant n^2}.$
Hence there is a
unique Leibniz deformation of
$\mathfrak{gl}(n).$
\end{example}

\begin{corollary}
Let $\mathcal{H}_N$ be the $(2N+1)$-dimensional complex Heisenberg Lie algebra
($N\geqslant 1$)
with basis $(x_i)_{1\leqslant i \leqslant 2N+1}$ and commutation relations
$[x_i,x_{N+i}] =x_{2N+1}$ $(1 \leqslant i \leqslant N).$
\\(i)
 $ZL^2_0(
{\mathcal{H}}_N,{\mathcal{H}}_N
 )$
 has basis
 $(x_{2N+1} \otimes (\omega^i \circ \omega^{j}))_{1 \leqslant i \leqslant j \leqslant 2N}$
 with
$(\omega_i)_{1\leqslant i \leqslant 2N+1}$ the dual basis to
$(x_i)_{1\leqslant i \leqslant 2N+1}.$
 \\(ii)
 $$\dim ZL^2_0({\mathcal{H}}_N , {\mathcal{H}}_N ) = \dim B^2(
{\mathcal{H}}_N,{\mathcal{H}}_N
 ) = N(2N+1);$$
 $$\dim HL^2(
{\mathcal{H}}_N,{\mathcal{H}}_N
 ) = \dim Z^2(
{\mathcal{H}}_N,{\mathcal{H}}_N
 ) =
 \begin{cases}
 \frac{N}{3}(8N^2 +6N+1) & \text{ if } N\geqslant 2 \\
 8 &\text{ if } N=1\, .
 \end{cases}
 $$
\end{corollary}
\begin{proof}
(i) Follows from
$\ker{\mathcal{I}} =
S^2  \left(
\mathfrak{g}/\mathcal{C}^2 \mathfrak{g}
\right)^*
.$
\\(ii)
First
$\mathcal{H}_N$
is adjoint $ZL^2$-uncoupling since it is $\mathcal{I}$-null.
The result then follows  from the fact that (\cite{commalg})
 $ \dim B^2(
{\mathcal{H}}_N,{\mathcal{H}}_N
 ) = N(2N+1)$ and for $N\geqslant 2,$
 $ \dim H^2(
{\mathcal{H}}_N,{\mathcal{H}}_N
 ) = \frac{2N}{3}(4N^2-1).$
\end{proof}

\begin{example}
\rm
The case $N=1$ has been studied in \cite{fialowski}.  In that case,
\linebreak[4]
 $\dim ZL^2_0({\mathcal{H}}_1 , {\mathcal{H}}_1 ) = 3$ and the 3 Leibniz deformations
 are nilpotent, in contradistinction with the 5 Lie deformations.



\end{example}

\end{document}